\documentclass[final,3p,11pt]{elsarticle}

\usepackage[utf8]{inputenc}
\usepackage{amsmath, amssymb, amsthm, bm, mathrsfs}
\usepackage{cases}
\usepackage{textcomp}
\usepackage{gensymb}
\usepackage[hidelinks]{hyperref}
\usepackage{array, tabularx, xcolor, colortbl}
\usepackage{csquotes, xpatch, algorithm, algpseudocode}
\usepackage{graphicx, color, float, subcaption, adjustbox}
\graphicspath{{Figures/}}
\usepackage{enumitem}
\usepackage[numbers]{natbib}
\definecolor{darkred}{rgb}{0.8,0,0}

\DeclareMathOperator{\vectorization}{\operatorname{vec}}

\begin{document}

\begin{frontmatter}
\title{On the fast assemblage of finite element matrices with application to nonlinear heat transfer problems}
\author[inst1]{Yannis Voet}
\ead{yannis.voet@epfl.ch}

\address[inst1]{MNS, Institute of Mathematics, École polytechnique fédérale de Lausanne, Station 8, CH-1015 Lausanne, Switzerland}

\date{\today}
\begin{abstract}
The finite element method is a well-established method for the numerical solution of partial differential equations (PDEs), both linear and nonlinear. However, the repeated reassemblage of finite element matrices for nonlinear PDEs is frequently pointed out as one of the bottlenecks in the computations. The second bottleneck being the large and numerous linear systems to be solved arising from the use of Newton's method to solve nonlinear systems of equations.
In this paper, we will address the first issue. We will see how under mild assumptions the assemblage procedure may be rewritten using a completely loop-free algorithm. Our approach leads to a small matrix-matrix multiplication for which we may count on highly optimized algorithms. 
\end{abstract}

\end{frontmatter}

\section{Introduction}
The mass and stiffness matrices are classical matrices arising when the finite element method (FEM) is used to discretize in space a PDE. Although the expression of the stiffness matrix differs from one application to another, it always depends on the gradients of the basis functions. Assembling these matrices may be computationally demanding, especially for nonlinear problems when they have to be reassembled several times. Hence, designing efficient matrix assemblage algorithms is critical. Several contributions have focused on designing algorithms of optimal or quasi-optimal complexity. An assemblage algorithm is optimal if $O(1)$ operations are required per nonzero entry. In \citep{ainsworth2011bernstein, ainsworth2016bernstein}, the authors designed an optimal complexity algorithm by choosing Bernstein polynomials as basis functions and making good use of their properties. Other researchers have exploited the tensor product structure of basis functions leading to sum factorization techniques \citep{antolin2015efficient} and weighted quadrature \cite{calabro2017fast, sangalli2018matrix}. Low-rank approximation techniques have also been proposed for tensor product basis functions in the context of isogeometric discretizations, allowing to approximate global matrices by a sum of Kronecker products and only assembling the (much smaller) factor matrices \citep{mantzaflaris2017low, scholz2018partial, hofreither2018black}. Finally, the geometry and physical nature of the problem may also be used to reduce the complexity in some circumstances \cite{antolin2020fast}. Besides complexity issues, another line of research has focused on designing algorithms that reconcile efficiency with the intrinsic ease of implementation of some programming languages. Indeed, classical assemblage algorithms require a loop over all finite elements. Although it may not be an issue for low level programming languages, it severely undermines the computational performance for interpreted languages such as MATLAB or Python. Thus, several attempts have been made over the last years to eliminate the loop over the elements by taking advantage of element-wise operations commonly defined in interpreted languages. Notable contributions include \cite{cuvelier2013efficient}, where the element matrices are vectorized and the loop over the elements is eliminated by taking advantage of some of the redundancy in the computations. However, their implementation uses a predefined quadrature formula leading to ad hoc order-specific algorithms. This lack of flexibility may be a shortcoming for nonlinear problems involving non-constant coefficients and requiring possibly different quadrature rules. Some of these issues were later addressed in \cite{cuvelier2016efficient} and included other interesting contributions to reduce the memory consumption of the algorithms and exploit the symmetry of the coefficient matrices being assembled. Although potential extensions to higher order discretizations were briefly discussed, the implementation remains quite technical and mostly left as future work. In \cite{rahman2013fast}, a tensorized approach is proposed (referred to as matrix-array operations) for the computation of local matrices for standard Lagrange finite elements. Their strategy was later extended in \cite{anjam2015fast} for Raviart-Thomas and Nédélec edge elements. The idea is again to compute all local matrices at once. However, the element matrices are not vectorized and instead local quantities are arranged along the pages of third order tensors. Specialized routines were coded up for standard linear algebra operations on the element level. These routines can now be replaced by a call to MATLAB's built-in \textit{pagemtimes} function introduced in the R2020b release. Nevertheless, we did not find any discussion related to nonlinear problems and non-constant coefficients. Similarly to the previous contributions, our strategy relies on the simultaneous computation of all local quantities. However, our approach is far more general as it supports any finite element order, any quadrature formula and non-constant coefficients while still being highly performing. The algorithms were designed from the very start with these requirements in mind, making them well-suited for nonlinear problems. For the sake of the exposition, we will restrict ourselves to scalar PDEs in a two dimensional setup with triangular meshes. However, the concepts are general and naturally carry over to a three dimensional setup. 

After recalling in Section \ref{se: fem_basics} the basis of the finite element method from which the matrices stem, the classical assemblage procedure is reviewed in Section \ref{se: classical_assemblage}. We then reformulate it in Section \ref{se: alternative_assemblage} for the mass and stiffness matrices. In Section \ref{se: nonlinear_fem}, we show how our formulation elegantly extends to tensors arising in nonlinear FEM methods. Section \ref{se: experiments} provides a few numerical experiments where the performance of our algorithm is compared to the classical assemblage algorithm and FreeFEM++. We conclude the experimental section with an application of FEM to nonlinear heat transfer problems where we compare the performance of our implementation to MATLAB's PDE toolbox. Potential extensions and future work is discussed in Section \ref{se: future_work}. Finally, Section \ref{se: conclusion} concludes with a summary of our findings.

\section{The classical finite element method} \label{se: fem_basics}
As a motivating example, we will consider a scalar parabolic PDE used for modeling diffusion-reaction processes. Let $\Omega \subset \mathbb{R}^d \ d=2,3$ be an open and bounded domain with Lipschitz boundary (such that the outward normal vector is well defined everywhere on the boundary). Let $I=[0, T]$ be the time domain with $T>0$ denoting the final time. Let us consider the following differential problem in strong form: find $u \colon \Omega \times [0, T] \to \mathbb{R}$ such that
\begin{align}
     m \frac{\partial u}{\partial t}-\nabla \cdot (c \nabla u) + a u &=f & &\text{ in } \Omega \times (0, \ T] \label{eq: diff_eq}\\
     u&=g_D & &\text{ on } \partial \Omega_D \times (0, \ T] \label{eq: dirichlet_bc}\\
     c \nabla u \cdot \mathbf{n} &= g_N & &\text{ on } \partial \Omega_N \times (0, \ T] \label{eq: neumann_bc}\\
     u(\mathbf{x},0)&=u_0(\mathbf{x}) & &\text{ in } \Omega \label{eq: initial_condition}
\end{align}
where $\partial \Omega_D$ and $\partial \Omega_N$ form a partition of the boundary $\partial \Omega$ such that $\overline{\partial \Omega_D} \cup \overline{\partial \Omega_N}=\partial \Omega$ and $\partial \Omega_D \cap \partial \Omega_N = \emptyset$. Equation \eqref{eq: diff_eq} is the differential equation to be solved. Equations \eqref{eq: dirichlet_bc} and \eqref{eq: neumann_bc} prescribe Dirichlet and Neumann boundary conditions, respectively, and Equation \eqref{eq: initial_condition} is an initial condition. We will assume the scalar coefficients $m$, $c$ and $a$ and the function $g_N$ may depend on the space and time variables as well as the unknown whereas the right-hand side function $f$ and the function $g_D$ depend only on the space and time variables. This model problem encompasses many engineering applications such as heat transfer problems and groundwater flow. We will make all the necessary assumptions on the data such that the problem stated is well-posed. We refer to classical textbooks such as \cite{quarteroni2009numerical}, Chapter 5 for a more detailed description.

The first step is to derive the weak form of the differential problem. We begin by multiplying Equation \eqref{eq: diff_eq} with a test function $v$. After using Green's identity and the divergence theorem, we obtain
\begin{equation}
    \int_{\Omega} m \frac{\partial u}{\partial t} v \ \mathrm{d} \Omega + \int_{\Omega} c \nabla u \cdot \nabla v \ \mathrm{d} \Omega + \int_{\Omega} a uv \ \mathrm{d} \Omega = \int_{\Omega} fv \ \mathrm{d} \Omega + \int_{\partial \Omega_N} g_N v \ \mathrm{d} S \quad \forall v \in V_0 \label{eq: weak_form}
\end{equation}
where we have defined $V_0=H_0^1=\{v \in H^1 \colon v|_{\partial \Omega_D} =0 \}$ as the space of test functions. Further defining the symmetric bilinear form
\begin{equation*}
    a(u,v)=\int_{\Omega} c \nabla u \cdot \nabla v \ \mathrm{d} \Omega + \int_{\Omega} a uv \ \mathrm{d} \Omega,
\end{equation*}
the problem now reads: find $u \colon \Omega \times [0, T] \to \mathbb{R}$ with $u(.,t) \in V=\{v \in H^1 \colon v|_{\partial \Omega_D} = g_D \}$ such that
\begin{equation*}
    \int_{\Omega} m \frac{\partial u}{\partial t} v \ \mathrm{d} \Omega + a(u,v) = \int_{\Omega} fv \ \mathrm{d} \Omega + \int_{\partial \Omega_N} g_N v \ \mathrm{d} S  \quad \forall v \in V_0.
\end{equation*}
The Galerkin method consists in searching for an approximation $u_h$ in a suitable finite dimensional space $V_h \subset V$. To ease the exposition, let us assume that $g_D=0$ such that $V=V_0$. Hence, we search for $u_h \colon \Omega \times [0, T] \to \mathbb{R}$ with $u_h(.,t) \in V_h$ such that
\begin{equation}
    \int_{\Omega} m \frac{\partial u_h}{\partial t} v_h \ \mathrm{d} \Omega + a(u_h,v_h) = \int_{\Omega} fv_h \ \mathrm{d} \Omega + \int_{\partial \Omega_N} g_N v_h \ \mathrm{d} S  \quad \forall v_h \in V_h. \label{eq: Galerkin_approx}
\end{equation}
In the finite element method, the approximation space $V_h$ is taken as the space of piecewise polynomials. Denoting $n=\dim(V_h)$ and $\{\phi_i\}_{i=1}^n$ the Lagrange basis for $V_h$, the discrete solution $u_h$ can be expressed as $u_h(\mathbf{x},t)=\sum_{i=1}^n u_i(t)\phi_i(\mathbf{x})$ where $\{u_i(t)\}_{i=1}^n$ are unknown time-dependent coefficients. It is enough to enforce Equation \eqref{eq: Galerkin_approx} on all basis functions $\{\phi_i\}_{i=1}^n$ since any $v_h \in V_h$ can be expressed as a linear combination of the basis functions. Hence, Equation \eqref{eq: Galerkin_approx} is equivalent to the set of equations
\begin{equation}
    \sum_{j=1}^n \frac{\partial u_j}{\partial t}(t) \int_{\Omega} m \ \phi_j \phi_i \ \mathrm{d} \Omega + \sum_{j=1}^n u_j(t) a(\phi_j,\phi_i) = \int_{\Omega} f \phi_i \ \mathrm{d} \Omega + \int_{\partial \Omega_N} g_N \phi_i \ \mathrm{d} S  \quad \forall i=1,\dots,n. \label{eq: indices_format}
\end{equation}
We now define the mass matrix $M_{ij}=\int_{\Omega} m \ \phi_j \phi_i \ \mathrm{d} \Omega$, $i,j=1,\dots,n$, the stiffness matrix $K_{ij}=a(\phi_j,\phi_i)$, $i,j=1\dots,n$, the right-hand side vector $\mathbf{f}(t)$ such that $f_i(t)=\int_{\Omega} f \phi_i \ \mathrm{d} \Omega + \int_{\partial \Omega_N} g_N \phi_i \ \mathrm{d} S$, $i=1,\dots,n$ and $\mathbf{u}(t)$ the vector of coefficients $\{u_i(t)\}_{i=1}^n$. Equation \eqref{eq: indices_format} can now be conveniently rewritten using vector notation as
\begin{equation}
    M\mathbf{\dot{u}}(t) + K\mathbf{u}(t) = \mathbf{f}(t) \label{eq: vector_format}
\end{equation}
where $M, K \in \mathbb{R}^{n \times n}$ and $\mathbf{f} \in \mathbb{R}^{n}$. Note that the stiffness matrix here accounts for two terms:
\begin{equation*}
    K_{ij}=a(\phi_j,\phi_i)=\int_{\Omega} c \ \nabla \phi_j \cdot \nabla \phi_i \ \mathrm{d} \Omega + \int_{\Omega} a \ \phi_j \phi_i \ \mathrm{d} \Omega.
\end{equation*}
Let us denote $C_{ij}=\int_{\Omega} c \  \nabla \phi_j \cdot \nabla \phi_i \ \mathrm{d} \Omega$ and $A_{ij}=\int_{\Omega} a \ \phi_j \phi_i \ \mathrm{d} \Omega$. From its definition, the matrix $A$ has the same structure as the matrix $M$. Hence, we will focus on the assemblage of the mass matrix $M$ and the conductivity matrix $C$. In the finite element literature, the matrix $C$ is also commonly called a stiffness matrix. However, the more general setup we have considered here required an extended terminology.

\section{The classical assemblage procedure} \label{se: classical_assemblage}
Part of the assemblage strategy is common to all matrices. Thus, let us consider the mass matrix for simplicity.
\begin{enumerate}
\item The assemblage procedure first takes advantage of the additivity of the integral such that
\begin{equation}
    M_{ij}=\int_{\Omega} m \ \phi_j \phi_i \ \mathrm{d} \Omega= \sum_{e=1}^{n_e} \int_{T_e} m \ \phi_j \phi_i \ \mathrm{d} \Omega \label{eq: assemblage_1}
\end{equation}
where $T_e$ for $e=1,\dots, n_e$ are the finite elements making up the mesh and $n_e$ is the number of elements. We will assume the domain is polyhedral such that it can be represented exactly by a union of triangles in 2D (tetrahedra in 3D).

\item The second step is to take advantage of the compact support of the basis functions. Indeed, the support of $\phi_i$ is the union of all triangles $T_k$ to which node $i$ belongs. Hence, 
$\int_{T_e} m \ \phi_j \phi_i \ \mathrm{d} \Omega=0$ if $i$ and $j$ do not belong simultaneously to the finite element $T_e$. Therefore, many of the terms (if not all) of the sum in Equation \eqref{eq: assemblage_1} are in fact zero for a given pair of indices $(i,j)$. It is this property that leads to sparse matrices. Hence, it is convenient to consider only the nonzero contributions obtained by restricting the indices $i$ and $j$ to those that belong simultaneously to $T_e$. Let $\mathcal{N}_e$ denote that set of indices. Therefore, for each triangle $T_e$, we compute the integrals
\begin{equation*}
    \int_{T_e} m \ \phi_j \phi_i \ \mathrm{d} \Omega \quad i,j \in \mathcal{N}_e \quad e=1,2,\dots,n_e.
\end{equation*}
The cardinality of $\mathcal{N}_e$ is equal to the number of nodes of the finite element. Let us denote it $n_p=|\mathcal{N}_e|$. Another way of viewing it is to define the local matrices $M_e \in \mathbb{R}^{n_p \times n_p}$ such that
\begin{equation}
    M_e=\int_{T_e} m \ \bm{\phi}_e \bm{\phi}_e^T \ \mathrm{d} \Omega \quad e=1,2,\dots,n_e \label{eq: local_mass_matrix}
\end{equation}
where the vector $\bm{\phi}_e$ contains the basis functions $\{\phi_i\}_{i \in \mathcal{N}_e}$.

\item Unfortunately, the domain of integration in Equation \eqref{eq: local_mass_matrix} is different for each finite element, which is not convenient for a computer implementation. Thus, we perform a change of variables. There are different options available. Here, we will consider an affine mapping between a reference triangle $\hat{T}$ and the current triangle $T_e$. This mapping is illustrated in Figure \ref{fig: affine_transformation}. We perform the change of variables $\mathbf{x}=\mathbf{F}_e(\mathbf{\hat{x}})=\mathbf{a}+B_e \mathbf{\hat{x}}$ where $B_e \in \mathbb{R}^{d \times d}$ is defined as $B_e=[\mathbf{b}-\mathbf{a}, \ \mathbf{c}-\mathbf{a}]$. This mapping has the advantage that its Jacobian matrix is constant and trivially given by $J_{F_e}=B_e$ regardless of which finite element order is used.

\begin{figure}[H]
    \centering
    \includegraphics[scale=0.5]{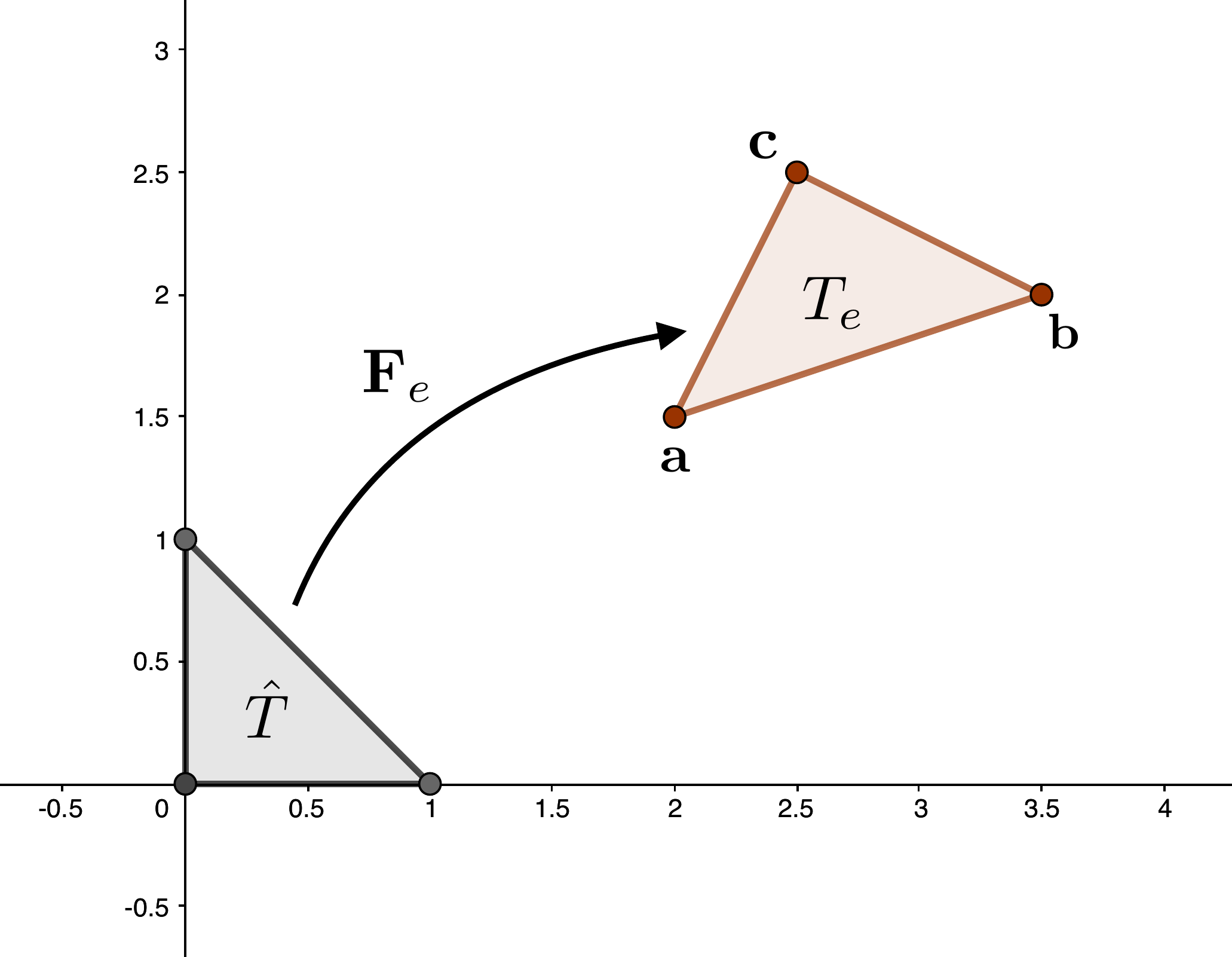}
    \caption{Affine transformation between the reference and the current triangle}
    \label{fig: affine_transformation}
\end{figure}

Denoting $\hat{m}=m(\mathbf{F}_e(\mathbf{\hat{x}}))$ and using the relation $\bm{\hat{\phi}}=\bm{\phi}_e(\mathbf{F}_e(\mathbf{\hat{x}}))$, we obtain
\begin{equation}
    M_e=\int_{\hat{T}} \hat{m} \ \bm{\hat{\phi}} \bm{\hat{\phi}}^T |\det(B_e)| \ \mathrm{d}\hat{\Omega} \quad e=1,2,\dots,n_e.
\end{equation}
The advantage is twofold. First of all, the integration domain no longer depends on the finite element. Secondly, the basis functions only need to be defined on the reference triangle. We refer to classical textbooks such as \cite{quarteroni2009numerical} for their expression depending on the finite element order.

\item Finally, integration is most efficiently done on a computer using numerical quadrature. In that respect, Gaussian quadrature which features the minimum number of quadrature nodes $n_q$ for a given degree of exactness is particularly attractive because it minimizes the number of function evaluations. Let $\{w_q\}_{q=1}^{n_q}$ and $\{\mathbf{\hat{x}}_q\}_{q=1}^{n_q}$ be the set of quadrature weights and nodes, respectively. Sets of quadrature weights and nodes are listed for example in \cite{laursen1978some} for different degrees of exactness. When numerical quadrature is used, the integral is approximated as follows

\begin{equation*}
    M_e \approx \sum_{q=1}^{n_q} w_q \hat{m}_q \ \bm{\hat{\phi}}_q \bm{\hat{\phi}}_q^T |\det(B_e)| \quad e=1,2,\dots,n_e
\end{equation*}
where we have denoted $\bm{\hat{\phi}}_q=\bm{\hat{\phi}}(\mathbf{\hat{x}}_q)$ and $\hat{m}_q=\hat{m}(\mathbf{\hat{x}}_q)$. 

The assemblage of the conductivity matrix $C$ is slightly more technical. Up to step 2, the procedure is the same and we consider
\begin{equation*}
    \int_{T_e} c \ \nabla \phi_j \cdot \nabla \phi_i \ \mathrm{d} \Omega \quad i,j \in \mathcal{N}_e \quad e=1,2,\dots,n_e.
\end{equation*}
We may rewrite it more compactly as
\begin{equation}
    C_e=\int_{T_e} c \ J_{\phi_e} J_{\phi_e}^T \ \mathrm{d} \Omega \quad e=1,2,\dots,n_e \label{eq: local_stiffness_matrix}
\end{equation}
where $J_{\phi_e}$ denotes the Jacobian matrix of the basis functions $\{\phi_i\}_{i \in \mathcal{N}_e}$.
At step 3, the change of variables leads to $J_{\phi_e}(\mathbf{F}_e(\mathbf{\hat{x}}))$. It must be related to the Jacobian matrix of the basis functions $\hat{\phi}$ defined on the reference triangle. For this purpose, let us recall the relation $\bm{\hat{\phi}}(\mathbf{\hat{x}})=\bm{\phi}_e(\mathbf{F}_e(\mathbf{\hat{x}}))$. Using the chain rule, we obtain
\begin{equation*}
    J_{\hat{\phi}}(\mathbf{\hat{x}})=J_{\phi_e}(\mathbf{F}_e(\mathbf{\hat{x}}))J_{F_e}(\mathbf{\hat{x}})=J_{\phi_e}(\mathbf{F}_e(\mathbf{\hat{x}}))B_e.
\end{equation*}
Therefore, we obtain the expression $J_{\phi_e}(\mathbf{F}_e(\mathbf{\hat{x}}))=J_{\hat{\phi}}(\mathbf{\hat{x}})B_e^{-1}$
and we note that $B_e$ is invertible provided the triangle is not degenerate.
After substituting in Equation \eqref{eq: local_stiffness_matrix}, we obtain
\begin{equation*}
    C_e=\int_{\hat{T}} \hat{c} \ J_{\hat{\phi}}(\mathbf{\hat{x}})(B_e^T B_e)^{-1} J_{\hat{\phi}}(\mathbf{\hat{x}})^T |\det(B_e)| \ \mathrm{d} \hat{\Omega} \quad e=1,2,\dots,n_e
\end{equation*}
where $\hat{c}=c(\mathbf{F}_e(\mathbf{\hat{x}}))$. Denoting $A_e=(B_e^T B_e)^{-1}$, using numerical quadrature and denoting $\hat{c}_q=\hat{c}(\mathbf{\hat{x}}_q)$, we end up with 
\begin{equation*}
    C_e \approx \sum_{q=1}^{n_q} w_q \hat{c}_q \ J_{\hat{\phi}}(\mathbf{\hat{x}}_q)A_e J_{\hat{\phi}}(\mathbf{\hat{x}}_q)^T |\det(B_e)| \quad e=1,2\dots,n_e.
\end{equation*}
\end{enumerate}
In the following, we shall drop the $\hat{}$ and denote $J_q=J_{\hat{\phi}}(\mathbf{\hat{x}}_q)$ for readability purposes. Thus, in a computer implementation, we set
\begin{align*}
    M_e &= \sum_{q=1}^{n_q} w_q m_q \ \bm{\phi}_q \bm{\phi}_q^T |\det(B_e)| & &\quad e=1,2,\dots,n_e, \\
    C_e &= \sum_{q=1}^{n_q} w_q c_q \ J_q A_e J_q^T |\det(B_e)| & &\quad e=1,2,\dots,n_e.
\end{align*}
The number of quadrature nodes $n_q$ usually depends on the matrix. Indeed, provided the coefficients $m$ and $c$ are constant, the expression of the mass matrix involves a higher degree polynomial compared to the stiffness matrix and one may want to use more quadrature nodes in order to integrate it exactly. Also note that thanks to the change of variables, $\bm{\phi}_q$ and $J_q$ only depend on the quadrature nodes and may be precomputed once and for all. Algorithm \ref{algo: classical_FEM_assemblage} is the classical assemblage algorithm. In practice, it is not exactly this version that is implemented but a variant which consists in first computing all the local matrices, storing them and then initializing the global matrices. In the process, the local matrices may be stored in a matrix of size $n_p \times n_p n_e$ for instance. This avoids having to either repeatedly add nonzero entries to sparse matrices or initializing a very large matrix full of zeros which might not fit in memory.

\begin{algorithm}[H]
\begin{algorithmic}[1]
\caption{Classical assemblage of FEM matrices}
\label{algo: classical_FEM_assemblage}
\Statex \textbf{Finite element parameters}:
\Statex $n_e$: number of finite elements
\Statex $n_p$: number of nodes per element
\Statex $n_q$: number of quadrature nodes per element
\Statex
\Statex \textbf{Input}: Quadrature weights $\{w_q\}_{q=1}^{n_q}$ and nodes $\{\mathbf{\hat{x}}_q\}_{q=1}^{n_q}$ potentially specific to each matrix
\Statex \textbf{Output}: Global mass matrix $M$ and conductivity matrix $C$
\Statex
\State Compute $\bm{\phi}_q$ and $J_q$ for $q=1,2,\dots,n_q$
\State Initialize $M$ and $C$
\For{$e=1,2,\dots,n_e$}
\State Initialize $M_e$ and $C_e$
\For{$q=1,2,\dots,n_q$}
\State $M_e \gets M_e+w_q m_q \ \bm{\phi}_q \bm{\phi}_q^T |\det(B_e)|$
\State $C_e \gets C_e + w_q c_q \ J_q A_e J_q^T |\det(B_e)|$
\EndFor
\State $M(\mathcal{N}_e, \mathcal{N}_e) \gets M(\mathcal{N}_e, \mathcal{N}_e) + M_e$
\State $C(\mathcal{N}_e, \mathcal{N}_e) \gets C(\mathcal{N}_e, \mathcal{N}_e) + C_e$
\EndFor
\State Return $M$ and $C$
\end{algorithmic}
\end{algorithm}  

\section{A reformulation of the assemblage procedure} \label{se: alternative_assemblage}
As we will see, our approach is nothing more than a reformulation of the same assemblage procedure which allows to eliminate both \textit{for} loops. The core strategy consists in not assembling the local matrices but their vectorization which allows to decouple the different dependencies and conveniently compute and store invariant quantities.

\subsection{Assemblage of the mass matrix}
\begin{figure}[H]
    \centering
    \includegraphics[scale=0.5]{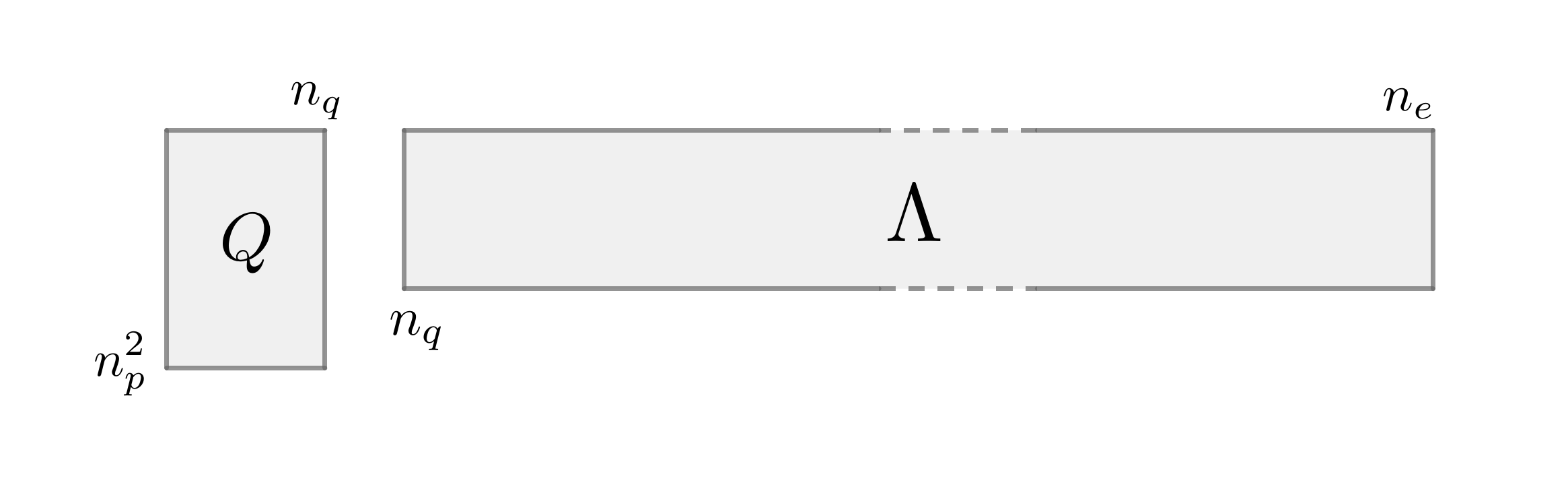}
    \caption{Illustration of the matrix sizes for the product $Q \Lambda$}
    \label{fig: QLambda_product}
\end{figure}

Let us first deal with the local mass matrix. Recall that it is given by
\begin{equation*}
    M_e = \sum_{q=1}^{n_q} w_q m_q \ \bm{\phi}_q \bm{\phi}_q^T |\det(B_e)|  \quad e=1,2,\dots,n_e.
\end{equation*}
Note that $w_q$, $m_q$ and $|\det(B_e)|$ are all scalars. Thus, let us regroup them and denote $\lambda_q=w_q m_q |\det(B_e)|$. Then, vectorizing the equation, we obtain
\begin{equation*}
    \vectorization(M_e)= \sum_{q=1}^{n_q} \lambda_q (\bm{\phi}_q  \otimes \bm{\phi}_q) \quad e=1,2,\dots,n_e.
\end{equation*}
Let us denote $\Phi=[\bm{\phi}_1, \bm{\phi}_2, \dots, \bm{\phi}_{n_q}] \in \mathbb{R}^{n_p \times n_q}$ and $\bm{\Lambda}_e=(\lambda_1, \lambda_2, \dots, \lambda_{n_q})^T \in \mathbb{R}^{n_q}$. Therefore, we obtain
\begin{equation*}
    \vectorization(M_e)=(\Phi \odot \Phi) \bm{\Lambda}_e \quad e=1,2,\dots,n_e
\end{equation*}
where $\odot$ denotes the Khatri-Rao product. For two matrices $A=[\mathbf{a}_1, \mathbf{a}_2, \dots, \mathbf{a}_m] \in \mathbb{R}^{n \times m}$ and $B=[\mathbf{b}_1, \mathbf{b}_2, \dots, \mathbf{b}_m]\in \mathbb{R}^{p \times m}$, their Khatri-Rao product is defined as
\begin{equation*}
    A \odot B = [\mathbf{a}_1 \otimes \mathbf{b}_1, \mathbf{a}_2 \otimes \mathbf{b}_2, \dots, \mathbf{a}_m \otimes \mathbf{b}_m] \in \mathbb{R}^{np \times m}.
\end{equation*}
Assuming all the elements of the mesh are of the same type, the vectorizations may be computed simultaneously since the matrix $(\Phi \odot \Phi)$ does not depend on the element. Thus, we obtain
\begin{equation*}
    [\vectorization(M_1), \vectorization(M_2), \dots, \vectorization(M_{n_e})]=(\Phi \odot \Phi) \Lambda = Q \Lambda
\end{equation*}
where we have denoted $Q=(\Phi \odot \Phi) \in \mathbb{R}^{n_p^2 \times n_q}$ and $\Lambda=[\bm{\Lambda}_1, \bm{\Lambda}_2, \dots, \bm{\Lambda}_{n_e}] \in \mathbb{R}^{n_q \times n_e}$. Therefore, the entire assemblage procedure is compactly written as a matrix-matrix multiplication between a small matrix $Q$ and a long matrix $\Lambda$ with only a few rows but many columns. From a linear algebra point of view, the matrix $Q$ can be seen as a basis for the vectorizations of the element matrices. Figure \ref{fig: QLambda_product} provides an illustration of the situation. The small matrix $Q$ can be easily computed and stored once and for all.

Let us now see how $\Lambda$ may be computed. We had denoted $\lambda_q=w_q m_q |\det(B_e)|$. Note that if all elements are of the same type, then the quadrature weights do not depend on the element. Moreover, if an affine mapping is used, $|\det(B_e)|$ does not depend on the quadrature nodes. Only the coefficients $m_q$ are both quadrature and element dependent. Making this dependency explicit, we have $\lambda_{qe}= m_{qe} w_q |\det(B_e)|$. Defining the vectors $\mathbf{w} \in \mathbb{R}^{n_q}$ and $\mathbf{b} \in \mathbb{R}^{n_e}$ such that $(\mathbf{w})_q=w_q$ and $(\mathbf{b})_e=|\det(B_e)|$ and the matrix $S_m \in \mathbb{R}^{n_q \times n_e}$ containing all the coefficient evaluations, then $\Lambda$ may be expressed as
\begin{equation}
    \Lambda= S_m \ast(\mathbf{w} \mathbf{b}^T) \label{eq: def_lambda}
\end{equation}
where $\ast$ denotes the Hadamard product. The Hadamard product of two matrices $A,B \in \mathbb{R}^{n \times m}$ is their element-wise product, defined as $(A \ast B)_{ij}=a_{ij}b_{ij}$. In practice, the computation of $\Lambda$ does not require forming the matrix $\mathbf{w} \mathbf{b}^T$ explicitly and then computing the Hadamard product. Since $\Lambda_{ij}=w_i({S_m})_{ij}b_j$, the products $w_i({S_m})_{ij}$ can be computed first by overwriting $S_m$ and then computing $w_i({S_m})_{ij}b_j$. Although the number of finite elements $n_e$ may be extremely large in applications, the number of quadrature nodes $n_q$ is typically very small and therefore computing these matrices does not lead to storage issues.

\subsection{Assemblage of the conductivity matrix}
Let us now investigate how the procedure may be adapted to the assemblage of the conductivity matrix. Recall that the local conductivity matrix was expressed as
\begin{equation*}
    C_e = \sum_{q=1}^{n_q} w_q c_q \ J_q A_e J_q^T |\det(B_e)|  \quad e=1,2,\dots,n_e.  
\end{equation*}
Analogously to what we did previously, we first define $\lambda_q=w_q c_q |\det(B_e)|$. Then, taking the vectorization, we obtain
\begin{equation*}
    \vectorization(C_e) = \sum_{q=1}^{n_q} \lambda_q \vectorization(J_q A_e J_q^T) = \sum_{q=1}^{n_q} \lambda_q (J_q \otimes J_q)\vectorization(A_e) \quad e=1,2,\dots,n_e.  
\end{equation*}
Now, after defining the matrix $Q=[J_1 \otimes J_1, J_2 \otimes J_2, \dots J_{n_q} \otimes J_{n_q}] \in \mathbb{R}^{n_p^2 \times d^2 n_q}$ and the vector
$\bm{\Lambda}_e=(\lambda_1, \lambda_2, \dots, \lambda_{n_q})^T \in \mathbb{R}^{n_q}$, the equation may be rewritten as
\begin{equation*}
    \vectorization(C_e) = Q (\bm{\Lambda}_e \otimes \vectorization(A_e)) \ \quad e=1,2,\dots,n_e. 
\end{equation*}
Now defining the matrix $W=[\vectorization(A_1), \vectorization(A_2), \dots, \vectorization(A_{n_e})] \in \mathbb{R}^{d^2 \times n_e}$ and gathering all the equations for $e=1,2,\dots,n_e$ we find
\begin{equation*}
    [\vectorization(C_1), \vectorization(C_2), \dots, \vectorization(C_{n_e})]= Q(\Lambda \odot W).
\end{equation*}
Once again, we end up with a simple matrix-matrix multiplication. The procedure to assemble the matrix $\Lambda$ is the same as before where the matrix $S_m$ must be replaced by the matrix $S_c$ containing the evaluations of the coefficient $c$. The matrix $Q$ defined here is slightly larger than the one defined for the mass matrix. In fact, it contains $d^2$ times more columns. However, since $d=2,3$ is a small integer, the matrix $Q$ remains very small and can be stored without any problem. We summarize this alternative assemblage procedure in Algorithm \ref{algo: alternative_FEM_assemblage}. Note that if the coefficients $m$ and $c$ are constant Line 5 can be skipped and Lines 6 and 7 reduce to $\Lambda_m=m \mathbf{w} \mathbf{b}^T$ and $\Lambda_c=c \mathbf{w} \mathbf{b}^T$, respectively. Whereas the assemblage includes the formation of the global matrices (Line 12), the novelty of our approach only lies in the formation of the element matrices (Line 1 to 11). From a computational point of view, the distinction between formation of element matrices and matrix assembly may be critical and will be specified when needed. However, for the sake of the presentation, we will often simply use the word assemblage to describe both.

\begin{algorithm}[H]
\begin{algorithmic}[1]
\caption{Alternative assemblage of FEM matrices}
\label{algo: alternative_FEM_assemblage}
\Statex \textbf{Finite element parameters}:
\Statex $n_e$: number of finite elements
\Statex $n_p$: number of nodes per element
\Statex $n_q$: number of quadrature nodes per element
\Statex \textbf{Input}:
\Statex Vector of quadrature weights $\mathbf{w}$ (potentially specific to each matrix) \Comment{$\mathbf{w} \in \mathbb{R}^{n_q}$}
\Statex Basis functions $\bm{\phi}_q=\bm{\hat{\phi}}(\mathbf{\hat{x}}_q)$ and Jacobian matrices $J_q=J_{\hat{\phi}}(\mathbf{\hat{x}}_q)$ evaluated at the quadrature nodes $\mathbf{\hat{x}}_q$ for $q=1,\dots,n_q$. \Comment{$\bm{\phi}_q \in \mathbb{R}^{n_p}$ and $J_q \in \mathbb{R}^{n_p \times d}$}
\Statex \textbf{Output}: 
\Statex Global mass matrix $M$ and conductivity matrix $C$
\Statex
\State Set $\Phi=[\bm{\phi}_1, \bm{\phi}_2, \dots, \bm{\phi}_{n_q}]$ \Comment{\parbox[t]{.25\linewidth}{$\Phi \in \mathbb{R}^{n_p \times n_q}$}}
\State Set $Q_m=(\Phi \odot \Phi)$ \Comment{\parbox[t]{.25\linewidth}{$Q_m \in \mathbb{R}^{n_p^2 \times n_q}$}}
\State Set $Q_c=[J_1 \otimes J_1, J_2 \otimes J_2, \dots J_{n_q} \otimes J_{n_q}]$ \Comment{\parbox[t]{.25\linewidth}{$Q_c \in \mathbb{R}^{n_p^2 \times d^2n_q}$}}
\State Compute the vector $\mathbf{b}$ and matrix $W$ from mesh related data. \Comment{\parbox[t]{.25\linewidth}{$\mathbf{b} \in \mathbb{R}^{n_e}$, $W \in \mathbb{R}^{d^2 \times n_e}$}}
\State Compute the matrices $S_m$ and $S_c$ containing coefficient evaluations \Comment{\parbox[t]{.25\linewidth}{$S_m, S_c \in \mathbb{R}^{n_q \times n_e}$}} 
\State Set $\Lambda_m=S_m \ast(\mathbf{w} \mathbf{b}^T)$ 
\State Set $\Lambda_c=S_c \ast(\mathbf{w} \mathbf{b}^T)$  \Comment{\parbox[t]{.25\linewidth}{$\Lambda_m, \Lambda_c \in \mathbb{R}^{n_q \times n_e}$}} 
\State Set $X_m=\Lambda_m$ \Comment{\parbox[t]{.25\linewidth}{$X_m \in \mathbb{R}^{n_q \times n_e}$}}
\State Set $X_c=\Lambda_c \odot W$ \Comment{\parbox[t]{.25\linewidth}{$X_c \in \mathbb{R}^{d^2n_q \times n_e}$}}
\State Set $V_m=Q_mX_m$ 
\State Set $V_c=Q_cX_c$ \Comment{\parbox[t]{.25\linewidth}{$V_m, V_c \in \mathbb{R}^{n_p^2 \times n_e}$}}
\State Form the matrices $M$ and $C$ from $V_m$ and $V_c$ and the set of indices of nonzero entries
\State Return $M$ and $C$
\end{algorithmic}
\end{algorithm}

\section{Consequences for nonlinear FEM} \label{se: nonlinear_fem}
In case one of the coefficients $m$, $c$ or $a$ or the function $g_N$ depends on the unknown, the problem becomes nonlinear and different solvers must be used. To cover such an event, let us rewrite the semi-discrete approximation in Equation \eqref{eq: vector_format} making the dependencies explicit
\begin{equation}
    M(\mathbf{u}(t)) \mathbf{\dot{u}}(t) + K(\mathbf{u}(t))\mathbf{u}(t) = \mathbf{f}(\mathbf{u}(t), t). \label{eq: semi_discrete_approx}
\end{equation}
Different numerical integration schemes may be used to solve approximately the resulting nonlinear system of differential equations. For simplicity, let us consider the implicit Euler method. Denoting $N$ the number of sub-intervals in time and $\Delta t=\frac{T}{N}$ the step size (assumed constant), the method requires solving the following nonlinear system of equations
\begin{equation*}
     M(\mathbf{u}^{s+1})(\mathbf{u}^{s+1}-\mathbf{u}^s)+\Delta t K(\mathbf{u}^{s+1}) \mathbf{u}^{s+1}=\Delta t \mathbf{f}(\mathbf{u}^{s+1}, t^{s+1}) 
\end{equation*}
for the unknown $\mathbf{u}^{s+1}$ for $s=0,1,\dots,N-1$ starting from $\mathbf{u}^0$ with $\mathbf{u}^s \approx \mathbf{u}(t^s)$, $t^s=s \Delta t \quad s=0,1,\dots,N$. Defining $\mathbf{F}(\mathbf{u})=M(\mathbf{u})(\mathbf{u}-\mathbf{u}^s)+\Delta t K(\mathbf{u}) \mathbf{u}-\Delta t \mathbf{f}(\mathbf{u}, t^{s+1})$, we solve the nonlinear system using Newton's method for which the iterations are defined as
\begin{equation*}
    \mathbf{u}_{l+1} = \mathbf{u}_l-\frac{ \partial \mathbf{F}(\mathbf{u}_l)}{\partial \mathbf{u}}^{-1}\mathbf{F}(\mathbf{u}_l) \quad l=0,1,2,\dots
\end{equation*}
until a prescribed tolerance is met. The sole question is how to compute the Jacobian matrix $\frac{ \partial \mathbf{F}(\mathbf{u})}{\partial \mathbf{u}}$. For this purpose, let us write the system component-wise
\begin{equation*}
    F_i(\mathbf{u})= \sum_{j=1}^n M_{ij}(\mathbf{u})(u_j-u_j^s)+\Delta t \sum_{j=1}^n K_{ij}(\mathbf{u})u_j - \Delta t f_i(\mathbf{u}, t^{s+1}) \quad i=1,\dots,n.
\end{equation*}
We now compute $\frac{\partial F_i}{\partial u_k}$ for $i,k=1,\dots,n$ considering the different terms
\begin{equation*}
    \frac{\partial}{\partial u_k}\Big(\sum_{j=1}^n M_{ij}(\mathbf{u})(u_j-u_j^s) \Big)= \sum_{j=1}^n \frac{\partial M_{ij}(\mathbf{u})}{\partial u_k}(u_j-u_j^s) + M_{ik}(\mathbf{u})=\sum_{j=1}^n \mathcal{M}_{ijk}(u_j-u_j^s) + M_{ik}(\mathbf{u}).
\end{equation*}
Similarly, we obtain
\begin{align*}
    &\frac{\partial}{\partial u_k}\Big( \sum_{j=1}^n K_{ij}(\mathbf{u})u_j \Big)= \sum_{j=1}^n \frac{\partial K_{ij}(\mathbf{u})}{\partial u_k}u_j + K_{ik}(\mathbf{u})=\sum_{j=1}^n \mathcal{K}_{ijk}u_j + K_{ik}(\mathbf{u}), \\
    &\frac{\partial}{\partial u_k}\Big(f_i(\mathbf{u}, t^{s+1}) \Big)= B_{ik}(\mathbf{u})
\end{align*}
and we defined the third order tensors
\begin{equation*}
    \mathcal{M}_{ijk}=\frac{\partial M_{ij}(\mathbf{u})}{\partial u_k} \text{ and } \mathcal{K}_{ijk}=\frac{\partial K_{ij}(\mathbf{u})}{\partial u_k}.
\end{equation*}
Their dependency on $\mathbf{u}$ is understood. Before giving their explicit expressions, we note that the discrete solution $u_h$ can be expressed as $u_h(\mathbf{x},t)=\bm{\phi}(\mathbf{x})^T \mathbf{u}(t)$ where the vector $\bm{\phi}(\mathbf{x})$ contains all the basis functions. Then, the tensors are expressed as
\begin{align*}
    &\mathcal{M}_{ijk}= \frac{\partial}{\partial u_k}\Big(\int_{\Omega} m(\bm{\phi}^T \mathbf{u}) \phi_i \phi_j \ \mathrm{d} \Omega \Big )
    =\int_{\Omega} m'(\bm{\phi}^T \mathbf{u}) \phi_i \phi_j \phi_k \ \mathrm{d} \Omega, \\
    &\mathcal{K}_{ijk}=\frac{\partial}{\partial u_k}\Big( \int_{\Omega} c(\bm{\phi}^T \mathbf{u}) \nabla \phi_i \cdot \nabla \phi_j+ a(\bm{\phi}^T \mathbf{u}) \phi_i \phi_j \ \mathrm{d} \Omega\Big)
    =\int_{\Omega} c'(\bm{\phi}^T \mathbf{u}) \nabla \phi_i \cdot \nabla \phi_j \phi_k + a'(\bm{\phi}^T \mathbf{u}) \phi_i \phi_j \phi_k \ \mathrm{d} \Omega.
\end{align*}
Finally,
\begin{equation*}
    B_{ik}(\mathbf{u}) = \frac{\partial}{\partial u_k}\Big(\int_{\Omega} f \phi_i \ \mathrm{d} \Omega + \int_{\partial \Omega_N} g_N(\bm{\phi}^T \mathbf{u}) \phi_i \ \mathrm{d} S \Big) 
    =\int_{\partial \Omega_N} g_N'(\bm{\phi}^T \mathbf{u}) \phi_i \phi_k \ \mathrm{d} S.
\end{equation*}
Regrouping all the expressions, we obtain
\begin{align*}
    \frac{\partial \mathbf{F}(\mathbf{u})}{\partial \mathbf{u}}=\mathcal{M} \circ_2 (\mathbf{u}-\mathbf{u}^n) + M(\mathbf{u})+ \Delta t ( \mathcal{K} \circ_2 \mathbf{u} + K(\mathbf{u})) - \Delta t B(\mathbf{u})
\end{align*}
where $\circ_2$ denotes the multiplication of a tensor and a vector in the second mode. Note that $\mathcal{M}, \mathcal{K} \in \mathbb{R}^{n \times n \times n}$ and $M,K,B \in \mathbb{R}^{n \times n}$. The tensor $\mathcal{M}$ is symmetric in all modes, meaning that all indices can be swapped. However, the tensor $\mathcal{K}$ is symmetric only in the two first modes ($\mathcal{K}_{ijk}=\mathcal{K}_{jik}$) due to the first term in its expression. Because of the degree of symmetry of these tensors, the multiplication can be done equivalently in the first mode. Unfortunately, when the multiplication of $\mathcal{K}$ and a vector is performed (in the first or second mode), the resulting matrix is nonsymmetric. Therefore, the Jacobian matrix $\frac{\partial \mathbf{F}(\mathbf{u})}{\partial \mathbf{u}}$ is in general nonsymmetric unless $c'=0$.

Having to deal with third order tensors may seem very worrisome, knowing that $n$ may reach over a million. However, the tensors are extremely sparse due to the compact support of the basis functions. Let us see how these tensors may be efficiently computed using the same arguments as for the mass and stiffness matrices. More specifically, we will be focusing on the assemblage of the mass tensor and conductivity tensor, which only considers the first term of the stiffness tensor.

\subsection{Assemblage of the mass tensor}
For readability purposes, let us denote $m'(\bm{\phi}^T \mathbf{u})$ simply $m'$, making the dependency on $\mathbf{u}$ implicit. Then,
\begin{equation*}
    \mathcal{M}_{ijk}=\int_{\Omega} m' \phi_i \phi_j \phi_k \ \mathrm{d} \Omega.
\end{equation*}
The exact same procedure described in Sections \ref{se: classical_assemblage} and \ref{se: alternative_assemblage} may be applied to compute this quantity. The same arguments as above lead to the computation of the local mass tensor
\begin{equation*}
    \mathcal{M}_e = \sum_{q=1}^{n_q} w_q m'_q \ \bm{\phi}_q \circ \bm{\phi}_q \circ \bm{\phi}_q |\det(B_e)| \quad e=1,2,\dots,n_e
\end{equation*}
where $\circ$ denotes the outer product. Denoting once more $\lambda_q=w_q m'_q |\det(B_e)|$ and taking the vectorization, we obtain
\begin{equation*}
    \vectorization(\mathcal{M}_e)=\sum_{q=1}^{n_q} \lambda_d \ (\bm{\phi}_q \otimes \bm{\phi}_q \otimes \bm{\phi}_q) = (\Phi \odot \Phi \odot \Phi) \mathbf{\Lambda}_e \quad e=1,2,\dots,n_e.
\end{equation*}
Gathering all the equations for $e=1,2,\dots,n_e$, and denoting $Q=(\Phi \odot \Phi \odot \Phi) \in \mathbb{R}^{n_p^3 \times n_q}$, we get
\begin{equation*}
    [\vectorization(\mathcal{M}_1), \vectorization(\mathcal{M}_2), \dots, \vectorization(\mathcal{M}_{n_e})]=(\Phi \odot \Phi \odot \Phi) \Lambda = Q \Lambda.
\end{equation*}
We recover once more a matrix-matrix multiplication. The matrix $Q$ is larger than the one defined for the assemblage of the mass matrix due to the fact that we are handling the vectorization of tensors. More specifically, it contains $n_p$ times more lines. Even if $n_p^3$ might become quite large for high order finite elements in 3D, $n_q$ typically remains small and thus storing $Q$ remains feasible.

\subsection{Assemblage of the conductivity tensor}
Similarly to before, we denote $c'(\bm{\phi}^T \mathbf{u})$ simply $c'$. Then, the conductivity tensor is expressed as
\begin{equation*}
    \mathcal{C}_{ijk}=\int_{\Omega} c' \nabla \phi_i \cdot \nabla \phi_j \phi_k \ \mathrm{d} \Omega.
\end{equation*}
Using the compact support of the basis functions followed by the same change of variables we obtain the local conductivity tensor
\begin{equation*}
    \mathcal{C}_e=\int_{\hat{T}} \hat{c}' \ J_{\hat{\phi}}(\mathbf{\hat{x}})(B_e^T B_e)^{-1} J_{\hat{\phi}}(\mathbf{\hat{x}})^T \circ \bm{\hat{\phi}} |\det(B_e)| \ \mathrm{d} \hat{\Omega} \quad e=1,2,\dots,n_e.
\end{equation*}
After approximating it using Gaussian quadrature, dropping the $\hat{}$ for readability, and denoting  $A_e=(B_e^T B_e)^{-1}$, we set
\begin{equation*}
   \mathcal{C}_e= \sum_{q=1}^{n_q} w_q c_q' \ J_q A_e J_q^T \circ \bm{\phi}_q |\det(B_e)| \quad e=1,2,\dots,n_e. 
\end{equation*}
As usual, denoting $\lambda_q=w_q c_q'|\det(B_e)|$ and taking the vectorization, we obtain the expression
\begin{equation*}
    \vectorization(\mathcal{C}_e)=\sum_{q=1}^{n_q} \lambda_q (\bm{\phi}_q \otimes \vectorization(J_q A_e J_q^T)) = \sum_{q=1}^{n_q} \lambda_q (\bm{\phi}_q \otimes J_q \otimes J_q)\vectorization(A_e).
\end{equation*}
Now, defining the matrix $Q=[\bm{\phi}_1 \otimes J_1 \otimes J_1, \bm{\phi}_2 \otimes J_2 \otimes J_2, \dots, \bm{\phi}_{n_q} \otimes J_{n_q} \otimes J_{n_q}] \in \mathbb{R}^{n_p^3 \times d^2 n_q}$, the equation may be rewritten as
\begin{equation*}
    \vectorization(\mathcal{C}_e) = Q (\bm{\Lambda}_e \otimes \vectorization(A_e)) \ \quad e=1,2,\dots,n_e. 
\end{equation*}
Denoting once more 
$W=[\vectorization(A_1), \vectorization(A_2), \dots, \vectorization(A_{n_e})] \in \mathbb{R}^{d^2 \times n_e}$ and gathering all the equations for $e=1,2,\dots,n_e$ we end up with the familiar equation
\begin{equation*}
    [\vectorization(\mathcal{C}_1), \vectorization(\mathcal{C}_2), \dots, \vectorization(\mathcal{C}_{n_e})]= Q(\Lambda \odot W).
\end{equation*}
This equation has the same structure as for the conductivity matrix. Only the definitions of the various matrices must be adjusted.

As we have seen, the assemblage of the mass and conductivity tensors does not induce any additional difficulty. Regardless of the object we are assembling, we end up with a matrix-matrix multiplication $QX$ where $Q$ is a small matrix that is computed once and for all and stored. This highly standardized assemblage procedure means Algorithm \ref{algo: alternative_FEM_assemblage} can in fact be used for the assemblage of the mass and conductivity tensors provided the definitions of the matrices involved are adjusted. The full potential of our procedure clearly appears for nonlinear problems. Indeed, decoupling the various quantities leads to tremendous savings of floating point operations since invariant quantities are only computed once. Let us summarize our findings at this stage. 

\begin{enumerate}[noitemsep]
    \item The computation of the nonzero entries of the finite element matrices and tensors always requires a matrix-matrix product $QX$.
    \item For the mass matrix and tensor $X=\Lambda$.
    \item For the conductivity matrix and tensor $X=\Lambda \odot W$.
    \item The matrix $\Lambda$ is always expressed as $\Lambda = S \ast(\mathbf{w} \mathbf{b}^T)$.
\end{enumerate}
For the assemblage of a given object (matrix or tensor):
\begin{enumerate}[noitemsep]
    \item The matrix $Q$ and the vector $\mathbf{w}$ only depend on the quadrature nodes and weights.
    \item The matrix $W$ and the vector $\mathbf{b}$ only depend on the mesh.
    \item The matrix $S$ is the only matrix to potentially depend on the unknown $\mathbf{u}$.
\end{enumerate}
Based on these findings, the matrix $S$ is therefore the only matrix to be recomputed at each step of Newton's method. Consequently, only a few products have to be carried out at each iteration to compute the nonzero entries. Provided the mesh does not change at each iteration, the matrices $Q$ and $W$ and the vectors $\mathbf{w}$ and $\mathbf{b}$ are only computed once. The price to pay is to store a few matrices. However, these matrices do not grow prohibitively with the size of the problem. Even on our largest test cases, storage issues were never experienced.

\section{Numerical experiments} \label{se: experiments}
\subsection{Experiment 1}
In this first experiment\footnote{The code used in the experiments is freely accessible at the following address: \url{https://c4science.ch/diffusion/12302/}}, we measure the performance of our assemblage algorithms by comparing them to FreeFEM++ \cite{hecht2012new}, an open-source and highly optimized finite element computer software written in C++, a low-level programming language. FreeFEM++ is a well-established reference in the finite element community and was also used by the authors in \citep{cuvelier2013efficient, cuvelier2016efficient} for performance comparison. In contrast, the implementation of Algorithms 1 and 2 is carried out in MATLAB, an interpreted high-level programming language. The unit square $\Omega=(0,1)^2$ is discretized by 6 increasingly fine meshes, generated using Gmsh \cite{geuzaine2009gmsh}, and whose mesh sizes are reported in Table \ref{tab: mesh_matrix_sizes} along with the sizes of the associated finite element matrices. The same meshes are used consistently in MATLAB and FreeFEM++ such that the matrix sizes coincide.

The computation times are compared for both $\mathbb{P}_1$ and $\mathbb{P}_2$ discretizations. The assemblage of finite element matrices for $\mathbb{P}_2$ discretizations is more demanding both because of the greater number of nodes and number of quadrature nodes per element needed to compute the integrals. In our experiments, the number of quadrature nodes is always taken to be the minimum number such that the integrals can be computed exactly provided the coefficients are constant. All the experiments are done on MacOS using a $2.2$ GHz Dual-Core Intel Core i7 processor with 8 GB of RAM. Since MATLAB does not provide a built-in function for computing the Khatri-Rao product of two matrices, we used the implementation from Laurent Sorber available on the Central File Exchange \cite{sorber2010khatri_rao}. 

\begin{table}[H]
    \begin{center}
    \begin{tabular}{|m{2cm}|m{3cm}|m{3cm}|}
    \hline
     Mesh size & Matrix size for $\mathbb{P}_1$ & Matrix size for $\mathbb{P}_2$ \\
     \hline
     $0.1 \times 2^{-1}$ & 568 & 2189\\
     $0.1 \times 2^{-2}$ & 2 211 & 8 681\\
     $0.1 \times 2^{-3}$ & 8 554 & 33 893\\
     $0.1 \times 2^{-4}$ & 33 803 & 134 573\\
     $0.1 \times 2^{-5}$ & 135 930 & 542 437\\
     $0.1 \times 2^{-6}$ & 542 812 & 2 168 685\\
    \hline
    \end{tabular}
    \caption{Mesh and matrix sizes}
    \label{tab: mesh_matrix_sizes}
    \end{center}
\end{table}

The average computation times over 5 consecutive launches are reported in Figures \ref{fig: timesDMatrix_p1} to \ref{fig: timesCMatrix_p2}. For our implementation of Algorithms 1 and 2, it is now important to distinguish the assemblage of the global matrices from the formation of all element matrices. Indeed, during the course of the experiments, we noticed that the call to MATLAB's built-in \textit{sparse} function at the end of the assemblage algorithm could be a major computational bottleneck. Thus, our results are reported separately with and without the function call. Unfortunately, this distinction cannot be done for FreeFEM++ unless the source code is modified, a risky meddling that requires a deep software knowledge. Thus, the computation times reported for FreeFEM++ are only for the assemblage of global matrices. Nevertheless, we believe element matrices are assembled into global matrices much more efficiently in FreeFEM++ than they are in MATLAB. Therefore, while comparing the assemblage of global matrices in FreeFEM++ to the formation of element matrices in MATLAB is certainly unfair, it does provide a very useful upper bound on the expected performance of our algorithms. On the other hand, comparing the assemblage of global matrices in both FreeFEM++ and MATLAB provides a safe lower bound on the expected performance. In other words, it gives an idea of the worst case performance.

The results reported in Figures \ref{fig: timesDMatrix_p1} to \ref{fig: timesCMatrix_p2} indicate that all algorithms have roughly the same complexity and scale linearly with the size of the matrices. This could be expected since Algorithm 2 is nothing more than a reformulation of Algorithm 1. As expected, FreeFEM++ is significantly faster than the standard assemblage algorithm (Algorithm 1), implemented in MATLAB. As a matter of fact, FreeFEM++ assembles the global matrices much faster than the element matrices in Algorithm 1, a clear testimony of the poor performance of embedded \textit{for} loops in MATLAB. On the contrary, Algorithm 2 reserves this trend completely. It is far ahead of FreeFEM++ for the formation of element matrices and interestingly preserves a very good margin for the assemblage of global matrices of sufficiently large size. It means that Algorithm 2 in its current state, even with the inefficient \textit{sparse} function, is significantly faster than FreeFEM++. We recall that it is essentially based on Khatri-Rao products and matrix-matrix multiplications for which there exists highly efficient implementations.

\begin{center}
\begin{minipage}[t]{.45\linewidth}
\vspace{0pt}
\begin{figure}[H]
    \includegraphics[width=\textwidth]{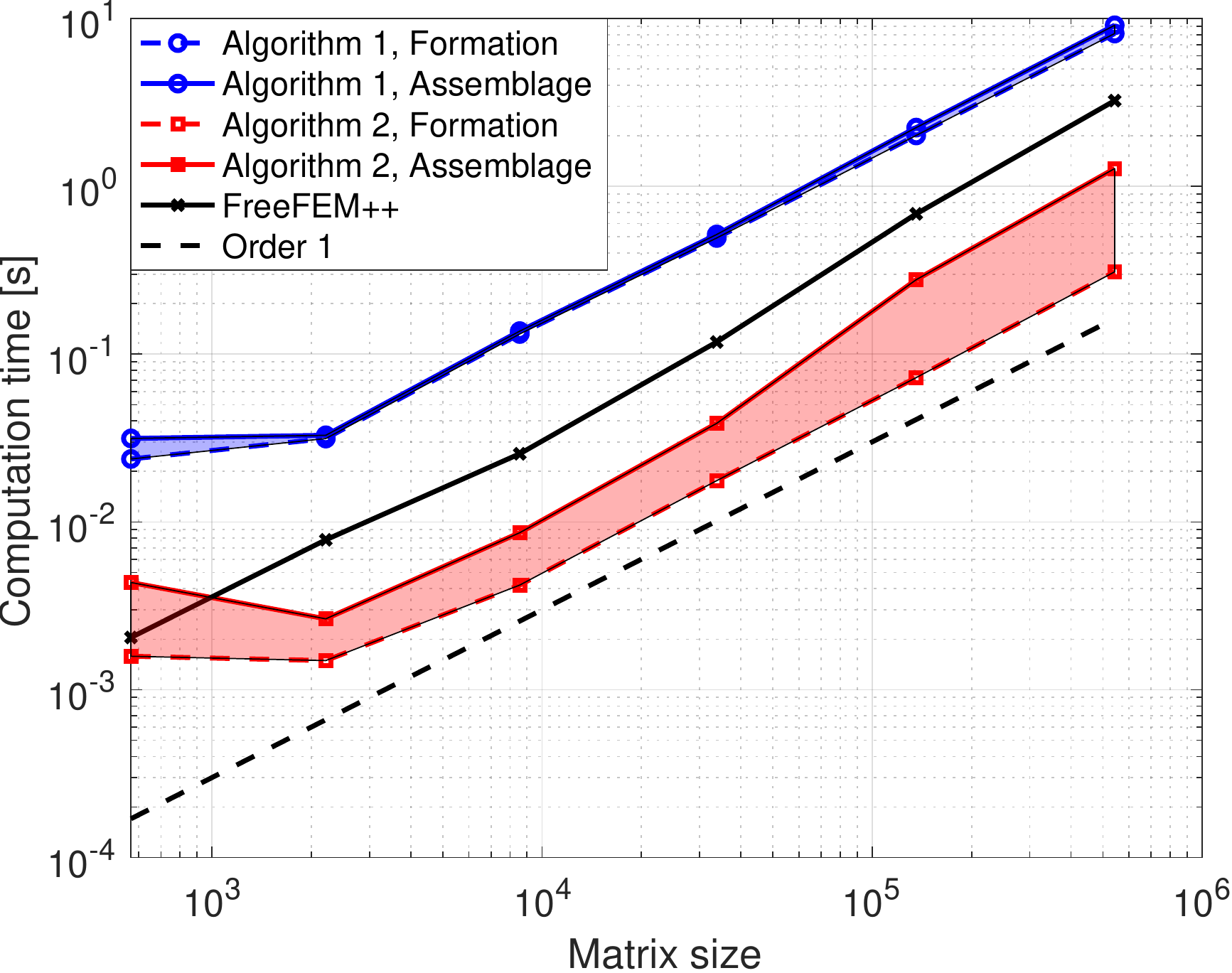}
    \caption{Computation times for the formation and assemblage of the mass matrix using $\mathbb{P}_1$ finite elements}
    \label{fig: timesDMatrix_p1}
\end{figure}
\end{minipage}
\hspace{2pt}
\begin{minipage}[t]{.45\linewidth}
\vspace{0pt}
\begin{figure}[H]
    \includegraphics[width=\textwidth]{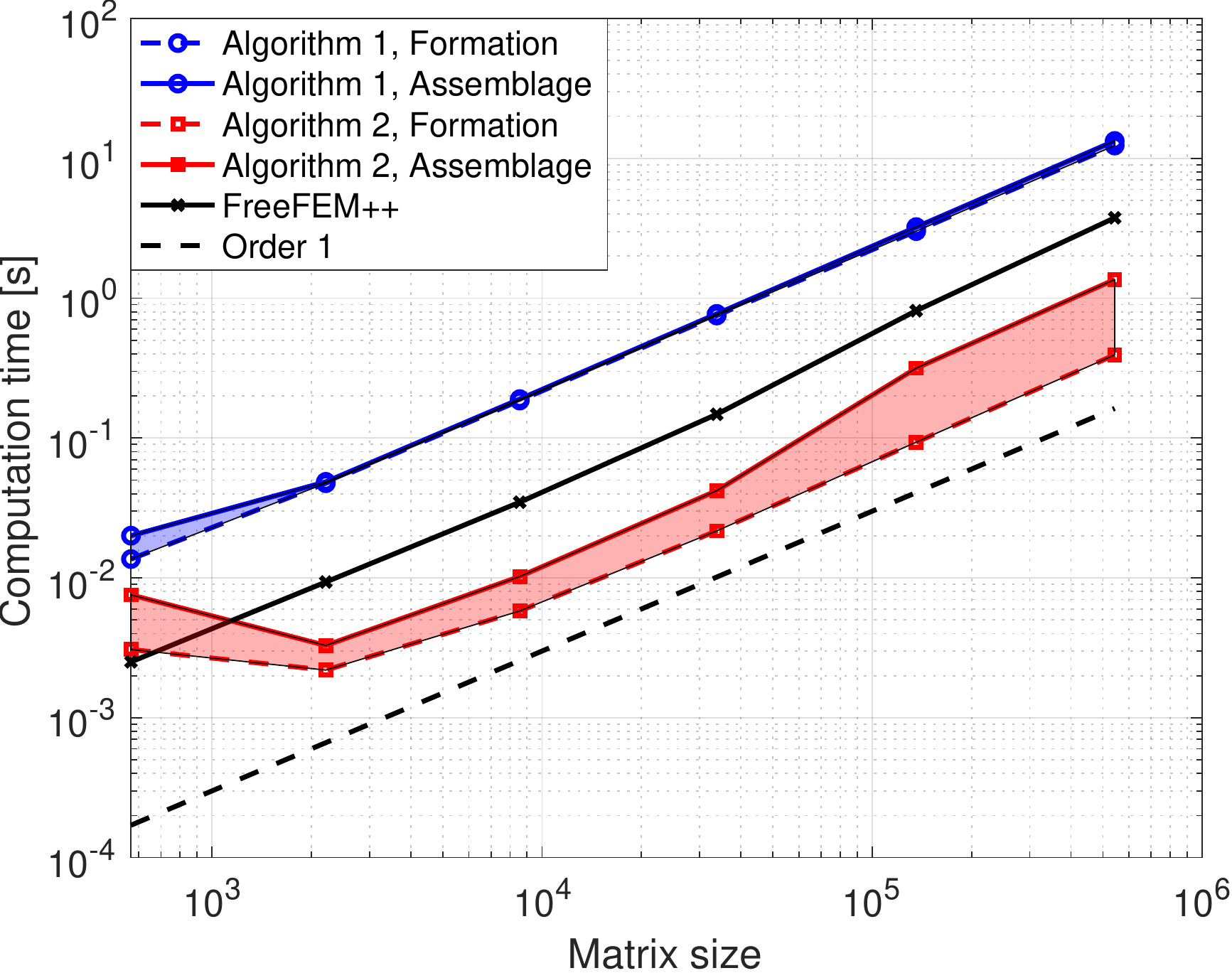}
    \caption{Computation times for the formation and assemblage of the conductivity matrix using $\mathbb{P}_1$ finite elements}
    \label{fig: timesCMatrix_p1}
\end{figure}
\end{minipage}
\end{center}

As a performance indicator, Table \ref{tab: factors_p1} reports the speedup factors for the formation of element matrices and the assemblage of global matrices for $\mathbb{P}_1$ discretizations. The results for $\mathbb{P}_2$ discretizations are reported in Table \ref{tab: factors_p2}. The speedup factor is defined as the ratio of computation times with respect to Algorithm 2. 
Since the computation time for the formation of element matrices cannot be easily retrieved in FreeFEM++, the ``Formation'' column actually divides the assemblage time of FreeFEM++ by the formation time of Algorithm 2 and provides an upper bound on the expected performance. The ``Assemblage'' column truly divides the assemblage times but since MATLAB's \textit{sparse} function takes such a large toll it provides a rather pessimistic (but conservative) lower bound on performance. According to this worst case analysis, Algorithm 2 is faster than FreeFEM++ by at least a factor 2.5 to 3.5 for $\mathbb{P}_1$ discretizations and 2.2 to 3.9 for $\mathbb{P}_2$ discretizations (and sufficiently large matrices). In general, we expect the speedup ratio to attain its maximum over a range of medium sized matrices. On the one hand, sufficiently large matrices are required for the computation time not to be dominated anymore by the least relevant parts of the algorithm such as the computation of the position of the nonzero entries. Since these parts of the algorithm are common to Algorithms 1 and 2, our strategy cannot make a significant difference. On the other hand, when the matrices become extremely large, MATLAB is forced to switch between different levels of cache memory. Large problems require a larger but also slower cache memory. Therefore, the speedup ratio slightly decreases but remains overall very large.

For Algorithm 1, profiling reports indicated that most of the computation time was spent in the inner \textit{for} loop over the quadrature nodes, as could be expected. For Algorithm 2, the bottleneck was found to be the call to MATLAB's built-in \textit{sparse} function at the end of the algorithm. The \textit{sparse} function's share in the overall computational expense in shown in Tables \ref{tab: time_frac_p1} and \ref{tab: time_frac_p2} for $\mathbb{P}_1$ and $\mathbb{P}_2$ discretizations, respectively. The toll taken by MATLAB's \textit{sparse} function is noticeable qualitatively by comparing the widths of the shaded areas in Figures \ref{fig: timesDMatrix_p1} to \ref{fig: timesCMatrix_p2}. The larger the width is, the greater the toll. It generally increases with the size of the matrix and may account for as much as 82\% of the computation time. For instance, for our largest benchmark with matrices of size over 2 million, all element mass matrices are computed in 1 s whereas the global mass matrix is assembled in 5.3 s.

\begin{table}[H]
\begin{center}
\begin{tabular}{l|c|c|c|c|c|c|c|c|}  
\arrayrulecolor{black} \cline{2-9}
 & \multicolumn{4}{c|}{Mass matrix} & \multicolumn{4}{c|}{Conductivity matrix} \\
\arrayrulecolor{black} \cline{2-9}
 & \multicolumn{2}{c|}{Formation} & \multicolumn{2}{c|}{Assemblage} & \multicolumn{2}{c|}{Formation} & \multicolumn{2}{c|}{Assemblage} \\
\hline
\multicolumn{1}{|l|}{Matrix size} & Algo 1 & FF++ & Algo 1 & FF++ & Algo 1 & FF++ & Algo 1 & FF++ \\
\hline
\multicolumn{1}{|l|}{568}    & 15.0    & 1.3  & 7.2 & 0.5 & 4.4 & 0.8 & 2.6 & 0.3    \\
\multicolumn{1}{|l|}{2 211}    & 21.0    & 5.2  & 12.4 & 2.9 & 21.7 & 4.3 & 14.8 & 2.8    \\
\multicolumn{1}{|l|}{8 554}    & 31.3    & 6.0  & 15.9 & 2.9 & 32.0  & 6.0 & 18.6 & 3.4    \\
\multicolumn{1}{|l|}{33 803}    & 27.9    & 6.7  & 13.3 & 3.0 & 34.7 & 6.8 & 18.4 & 3.5    \\
\multicolumn{1}{|l|}{135 930}    & 27.7    & 9.5  & 8.1 & 2.5 & 32.3 & 8.7 & 10.2 & 2.6    \\
\multicolumn{1}{|l|}{542 812}    & 26.3    & 10.5  & 7.1 & 2.5 & 31.4 & 9.5 & 9.8 & 2.8    \\
\hline
\end{tabular}
\caption{Speedup factors for Algorithm 1 (Algo 1) and FreeFEM++ (FF++) for the formation of element matrices and assemblage of global matrices using $\mathbb{P}_1$ finite elements.}
\label{tab: factors_p1}
\end{center}
\end{table}

\begin{table}[H]
\begin{center}
\begin{tabular}{l|c|c|c|c|}  
\cline{2-5}
 & \multicolumn{2}{c|}{Mass matrix} & \multicolumn{2}{c|}{Conductivity matrix} \\
\hline
\multicolumn{1}{|l|}{Matrix size} & Algorithm 1 & Algorithm 2 & Algorithm 1 & Algorithm 2 \\
\hline
\multicolumn{1}{|l|}{568}    & 24.5\%  & 63.6\%  & 31.9\% & 59.3\%  \\
\multicolumn{1}{|l|}{2 211}  & 4.4\%   & 43.6\%  & 2.2\% & 33.3\% \\
\multicolumn{1}{|l|}{8 554}   & 3.8\%    & 51.2\%  & 2.3\% & 43.3\% \\
\multicolumn{1}{|l|}{33 803}    & 4.7\%    & 54.5\%  & 2.6\% & 48.5\% \\
\multicolumn{1}{|l|}{135 930}    & 10.2\%    & 73.9\%  & 6.4\% & 70.4\% \\
\multicolumn{1}{|l|}{542 812}    & 10.6\%    & 75.7\%  & 7.3\% & 71.1\% \\
\hline
\end{tabular}
\caption{Time fraction spent in the \textit{sparse} function for the assemblage of global matrices using $\mathbb{P}_1$ finite elements}
\label{tab: time_frac_p1}
\end{center}
\end{table}

\begin{center}
\begin{minipage}[t]{.45\linewidth}
\vspace{0pt}
\begin{figure}[H]
    \includegraphics[width=\textwidth]{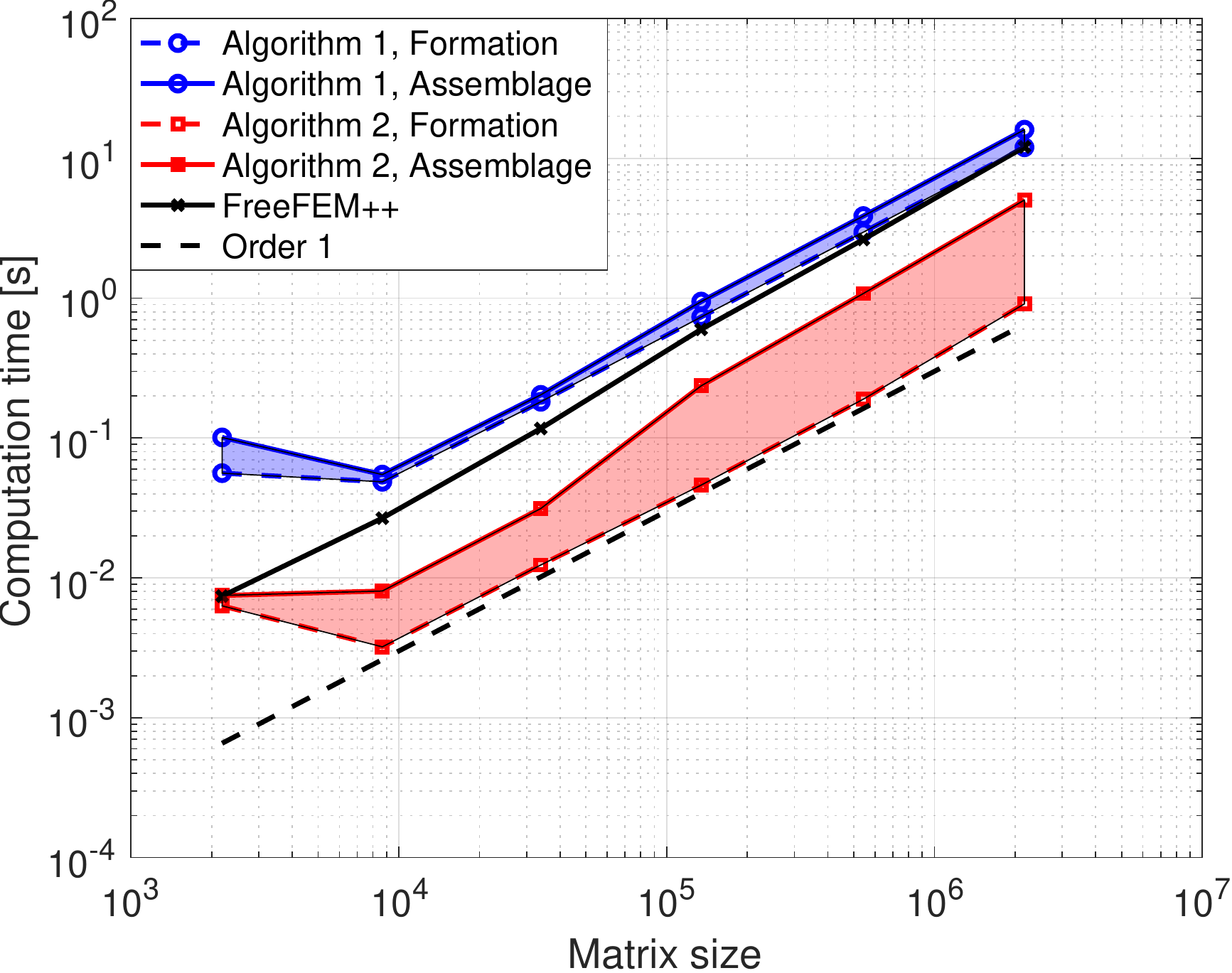}
    \caption{Computation times for the formation and assemblage of the mass matrix using $\mathbb{P}_2$ finite elements}
    \label{fig: timesDMatrix_p2}
\end{figure}
\end{minipage}
\hspace{2pt}
\begin{minipage}[t]{.45\linewidth}
\vspace{0pt}
\begin{figure}[H]
    \includegraphics[width=\textwidth]{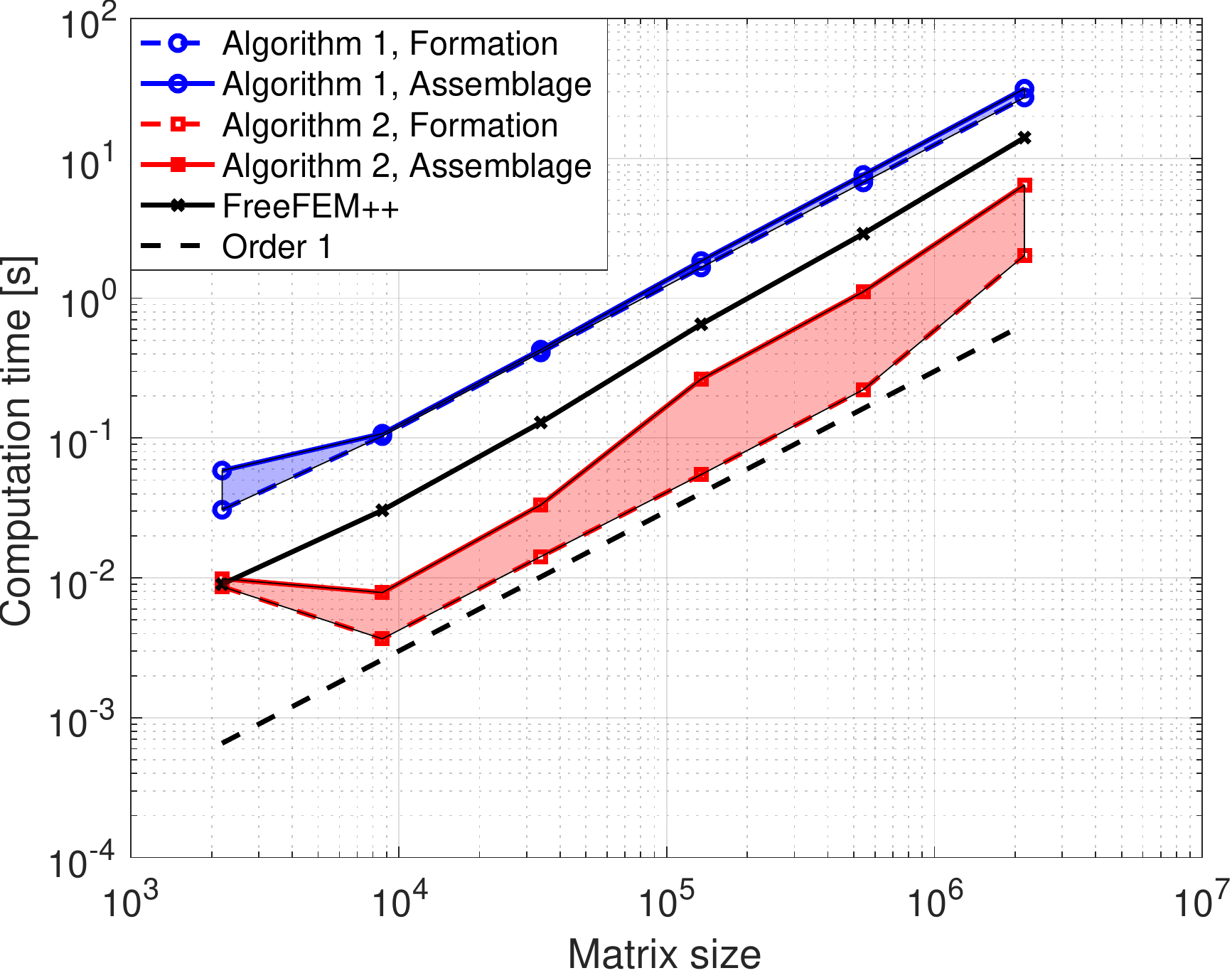}
    \caption{Computation times for the formation and assemblage of the conductivity matrix using $\mathbb{P}_2$ finite elements}
    \label{fig: timesCMatrix_p2}
\end{figure}
\end{minipage}
\end{center}

\begin{table}[H]
\begin{center}
\begin{tabular}{l|c|c|c|c|c|c|c|c|}  
\arrayrulecolor{black} \cline{2-9}
 & \multicolumn{4}{c|}{Mass matrix} & \multicolumn{4}{c|}{Conductivity matrix} \\
\arrayrulecolor{black} \cline{2-9}
 & \multicolumn{2}{c|}{Formation} & \multicolumn{2}{c|}{Assemblage} & \multicolumn{2}{c|}{Formation} & \multicolumn{2}{c|}{Assemblage} \\
\hline
\multicolumn{1}{|l|}{Matrix size} & Algo 1 & FF++ & Algo 1 & FF++ & Algo 1 & FF++ & Algo 1 & FF++ \\
\hline
\multicolumn{1}{|l|}{2 189}    & 8.9    & 1.2  & 13.5 & 1 & 3.5 & 1.0 & 5.9 & 0.9    \\
\multicolumn{1}{|l|}{8 681}    & 15.1    & 8.3  & 6.8 & 3.3 & 28.1 & 8.3 & 13.7 & 3.9    \\
\multicolumn{1}{|l|}{33 893}    & 14.7    & 9.5  & 6.5 & 3.7 & 28.8  & 9.1 & 12.9 & 3.9    \\
\multicolumn{1}{|l|}{134 573}    & 15.9    & 13.0  & 4.0 & 2.5 & 30.3 & 11.9 & 7.0 & 2.5    \\
\multicolumn{1}{|l|}{542 437}    & 15.6    & 13.9  & 3.6 & 2.4 & 30.3 & 13.0 & 6.8 & 2.6    \\
\multicolumn{1}{|l|}{2 168 685}    & 13.2    & 13.2  & 3.2 & 2.4 & 13.4 & 7.0 & 4.9 & 2.2    \\
\hline
\end{tabular}
\caption{Speedup factors for Algorithm 1 (Algo 1) and FreeFEM++ (FF++) for the formation of element matrices and assemblage of global matrices using $\mathbb{P}_2$ finite elements.}
\label{tab: factors_p2}
\end{center}
\end{table}

\begin{table}[H]
\begin{center}
\begin{tabular}{l|c|c|c|c|}  
\cline{2-5}
 & \multicolumn{2}{c|}{Mass matrix} & \multicolumn{2}{c|}{Conductivity matrix} \\
\hline
\multicolumn{1}{|l|}{Matrix size} & Algorithm 1 & Algorithm 2 & Algorithm 1 & Algorithm 2 \\
\hline
\multicolumn{1}{|l|}{2 189}    & 44.5\%  & 15.7\%  & 47.5\% & 12.3\%  \\
\multicolumn{1}{|l|}{8 681}  & 11.2\%   & 60.1\%  & 4.0\% & 53.3\% \\
\multicolumn{1}{|l|}{33 893}   & 10.8\%    & 60.5\%  & 4.4\% & 57.3\% \\
\multicolumn{1}{|l|}{134 573}    & 22.4\%    & 80.5\%  & 10.0\% & 79.2\% \\
\multicolumn{1}{|l|}{542 437}    & 23.3\%    & 82.4\%  & 11.6\% & 80.0\% \\
\multicolumn{1}{|l|}{2 168 685}    & 25.2\%    & 81.9\%  & 13.4\% & 68.6\% \\
\hline
\end{tabular}
\caption{Time fraction spent in the \textit{sparse} function for the assemblage of global matrices using $\mathbb{P}_2$ finite elements}
\label{tab: time_frac_p2}
\end{center}
\end{table}

\subsection{Experiment 2}
In this second experiment, we propose to test our implementation on a nonlinear transient heat transfer problem. Computationally speaking, nonlinear problems are very demanding both because of the repeated reassemblage of finite element matrices and because of the numerous linear systems to be solved. The case study is taken from \cite{ferreira2017validation}, where a series of validation examples for thermal finite element software are reported. The geometry is a square of side length $0.2$ m, $\Omega=(0,0.2)^2$. The source of non-linearity is coming from both a temperature dependent thermal conductivity and radiation boundary conditions. The material properties, geometry and boundary conditions are summarized in Table \ref{tab: material_example2} and Figure \ref{fig: geometry_example2}.

\begin{figure}[H]
    \centering
    \includegraphics[scale=0.47]{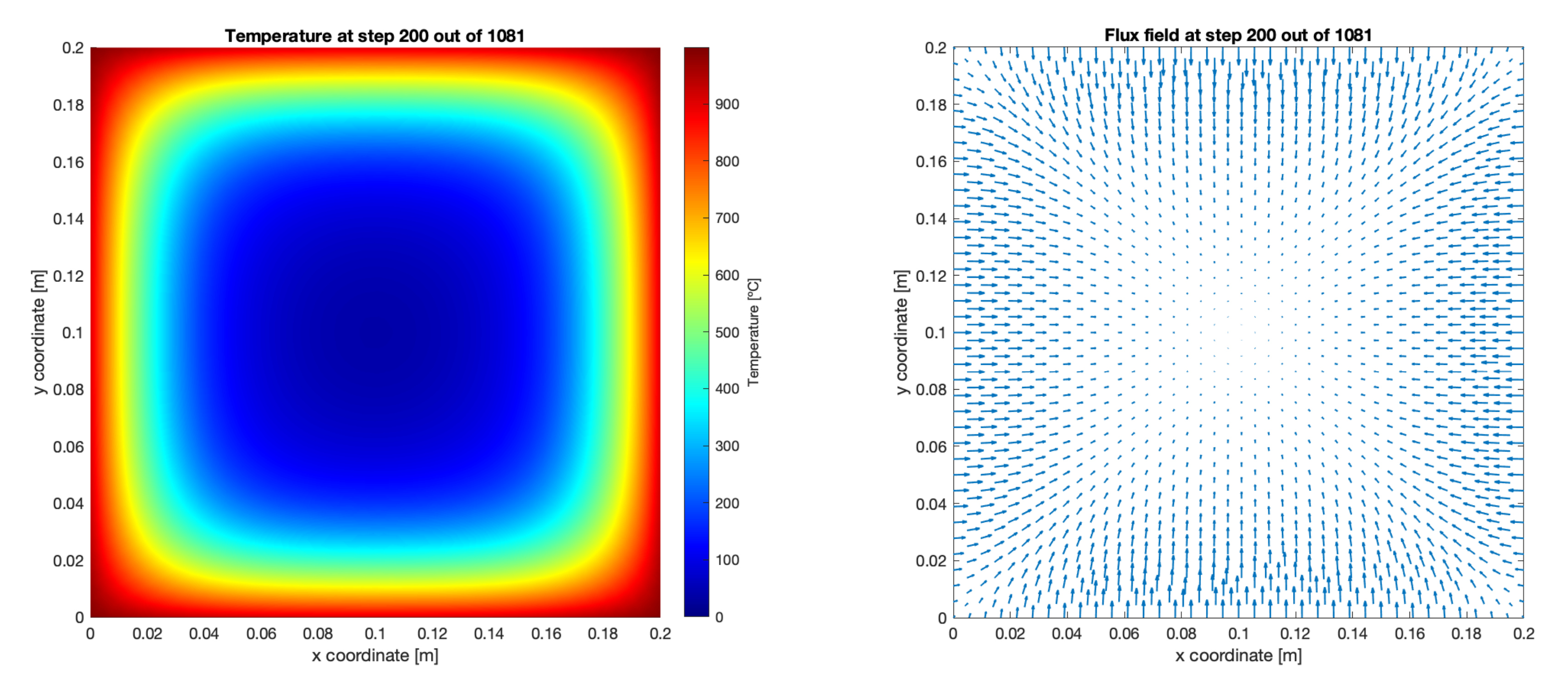}
    \caption{Finite element solution at time $t= 1990$ s}
    \label{fig: example2_all}
\end{figure}

\begin{figure}[H]
    \centering
    \includegraphics[scale=0.6]{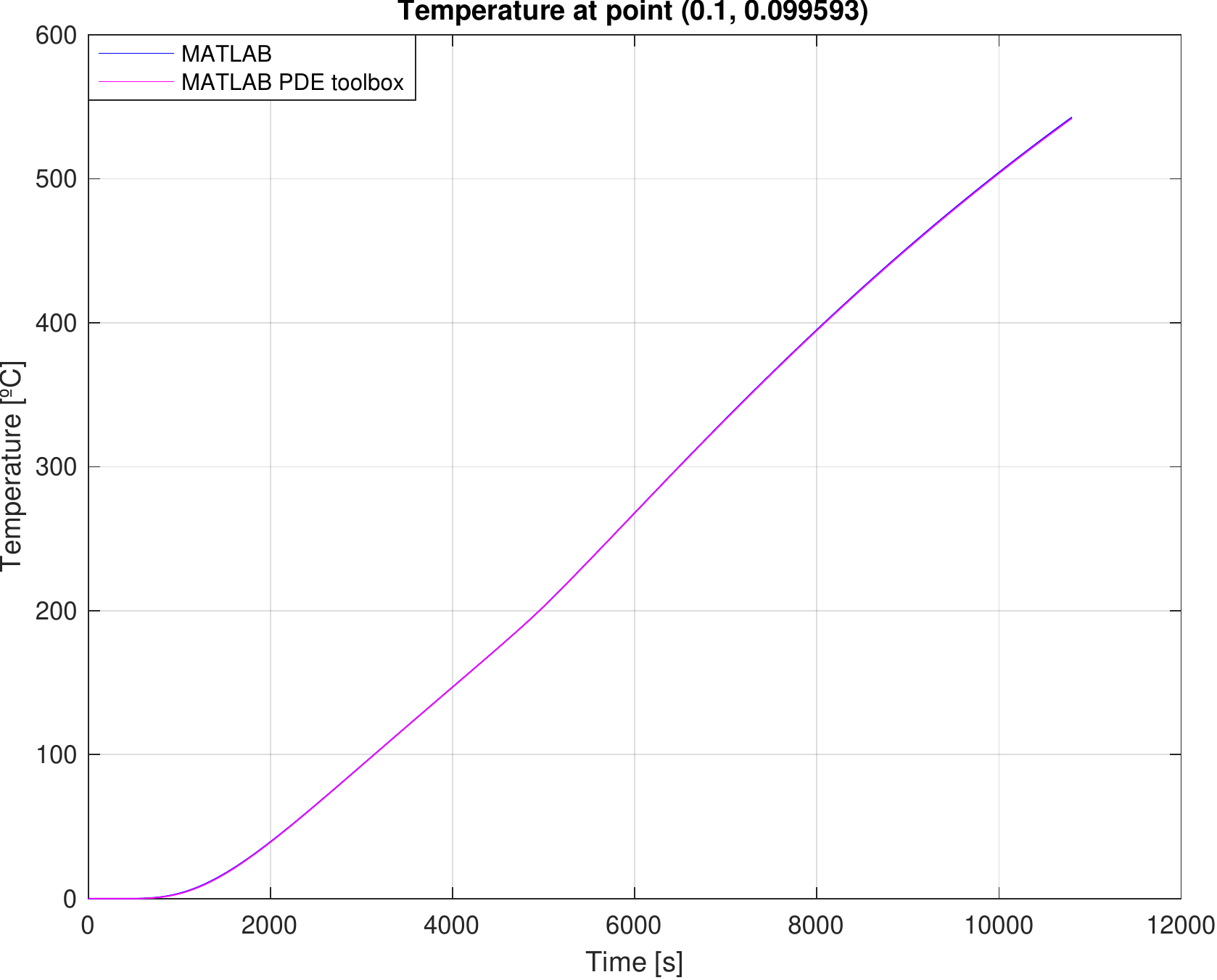}
    \caption{Temperature curves obtained using our in-house solver and MATLAB's PDE toolbox}
    \label{fig: example2_comparison}
\end{figure}

\begin{minipage}[t]{.2\linewidth}
\vspace{0pt}
\begin{table}[H]
    \begin{tabular}{|m{2.5cm}|m{1cm}|}
    \hline
     Quantity & Value \\
     \hline
     $c \ [\frac{\text{J}}{\text{kgK}}]$ & 1000 \\
     $\rho \ [\frac{\text{kg}}{\text{m\textsuperscript{3}}}]$ & 2400 \\
     thickness [m] & 1 \\
     height [m] & 0.2 \\
     length [m] & 0.2 \\
    \hline
    \end{tabular}
    \caption{Material and geometrical properties}
    \label{tab: material_example2}
\end{table}
\end{minipage}
\hspace{2cm}
\begin{minipage}[t]{.65\linewidth}
\vspace{0pt}
\begin{figure}[H]
    \centering
    \includegraphics[width=\textwidth]{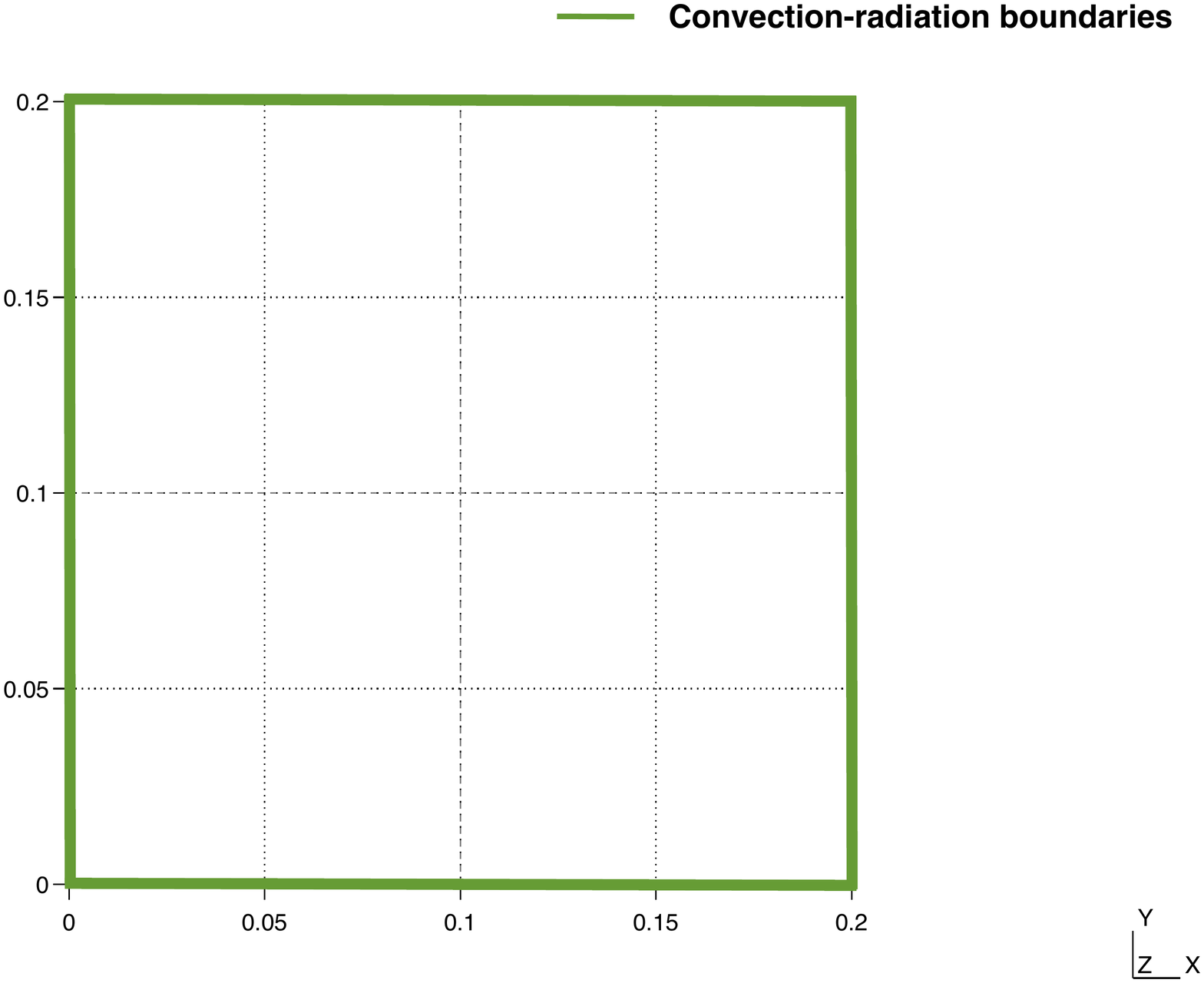}
    \caption{Geometry and boundary conditions}
    \label{fig: geometry_example2}
\end{figure}
\end{minipage}
The thermal conductivity is a piecewise linear function such that
\begin{equation*}
    k(T)=
\begin{cases}
k_1+\frac{k_2-k_1}{T_2-T_1}(T-T_1) & T_1 \leq T\leq T_2 \\
k_2+\frac{k_3-k_2}{T_3-T_2}(T-T_2) & T_2 < T\leq T_3
\end{cases}
\end{equation*}
with $T_1=0 \ \degree C$, $T_2=200 \ \degree C$, $T_3=1000 \ \degree C$ and $k_1=1.5$, $k_2=0.7$, $k_3=0.5$. All thermal conductivities are expressed in $\frac{\text{W}}{\text{mK}}$. The differential problem reads

\begin{align*}
     c \rho \frac{\partial T}{\partial t}-\nabla \cdot (k \nabla T) &=0 & &\text{ in } \Omega \times (0, \ T]\\
     k \nabla T \cdot \mathbf{n} &= h_c(T_a-T)+h_r(T_a^4-T^4) & &\text{ on } \partial \Omega \times (0, \ T]\\
     T(\mathbf{x}, 0)&=T_0 & &\text{ in } \Omega
\end{align*}
with coefficients $h_c=10 \ \frac{\text{W}}{\text{m\textsuperscript{2}K}}$ and $h_r=\epsilon \sigma$ with $\epsilon=0.8$ being the emissivity and $\sigma=5.670373\times 10^{-8} \ \frac{\text{W}}{\text{m\textsuperscript{2}K\textsuperscript{4}}}$ the Stefan-Boltzmann constant. The ambient temperature is set to $T_a=1000 \ \degree C$ and the initial temperature is $T_0 = 0 \ \degree C$. Finally, the simulation is carried out over three hours ($T=10800$ s). The time step was set to $\Delta t=10$ s. The geometry was discretized using $\mathbb{P}_2$ finite elements and the implicit Euler method was used for the numerical integration of the nonlinear system of differential equations. The nonlinear systems were solved using the classical Newton method described in Section \ref{se: nonlinear_fem}.
The simulations were run independently on our in-house solver and MATLAB's PDE toolbox. A similar mesh size was used in both cases such that the meshes generated contained 7593 and 7465 nodes, respectively. The finite element solution at time $t= 1990$ s along with the heat flux vectors is reported in Figure \ref{fig: example2_all} for illustration. To compare the results of our in-house solver with those of MATLAB's PDE toolbox, the temperature increase at the center of the domain as a function of time is reported in Figure \ref{fig: example2_comparison}. The curves overlap perfectly and validate our implementation.

To compare the performance of the implementations, we have generated profiling reports. Finite element software typically involve multiple stages from pre-processing to post-processing. In order to make a fair comparison, we have restricted our attention to the computation times needed to compute the finite element solution. More specifically, we refer to it as the time spent in the \textit{solve} function of the program. Although we do not know which scheme is being used in MATLAB both for the numerical integration and the solution to the nonlinear systems, we expect it to use state-of-the-art implementations. In comparison, our implementation uses a basic numerical integration scheme and the classical Newton method. The latter is known to be far too expensive and is not much used in practice. We intentionally decided to carry on with it because it provides a good test for our algorithm. On the other hand, we carefully tried to optimize the remaining part of our program. In particular, knowing that calls to MATLAB's \textit{sparse} function were the bottlenecks in Algorithm 2, we tried to use it as little as possible. For this purpose, we designed linear operators to carry out matrix-vector or tensor-vector multiplications without explicitly assembling the matrix nor the tensor. Secondly, we designed an algorithm to detect the various sources of non-linearity and hence avoid unnecessary computations. The idea being to only recompute quantities which are changing. We believe MATLAB uses a similar strategy. The highly standardized assembly procedure we have described allows us to assemble all necessary matrices in a single file. Algorithm 2 is the workhorse behind the assembly procedure. We ran our simulation twice using two different strategies to solve the linear systems of equations. 
\begin{enumerate}
    \item The linear systems are solved using the backslash (\textbackslash) command.
    
    \item The linear systems are solved using an in-house implementation of GMRES \cite{saad1986gmres}.
\end{enumerate}
The first strategy is obviously very bad since a matrix factorization will be recomputed at each step of the Newton iterations. Nevertheless, the Newton method converged very rapidly and required only $2$ to $3$ iterations for each nonlinear system to meet the prescribed tolerance on the norm of the increment at $10^{-7}$. Profiling reports indicated that in total $2757$ linear systems were solved for $1080$ time-steps. Despite the tremendous number of linear systems, our implementation outperformed MATLAB's PDE toolbox by a factor $2.07$, the computation times being $186.1$ s against $385.1$ s. Due to the poor choice of solver for the nonlinear and linear systems, it is reasonable to believe that the performance of our implementation comes from the reassemblage procedure. Profiling reports indicated that about $74\%$ of the computation time was spent solving the linear systems. The formation of element matrices only accounted for $7.6\%$ of the computation time.

In a second experiment, we chose a slightly better strategy to solve the linear systems. We used the iterative method GMRES. This choice stems from the fact that the linear systems involve a nonsymmetric matrix in general. Despite the large matrix size, its good conditioning allowed to solve the linear systems to the prescribed tolerance of $10^{-7}$ within $40$ iterations without any preconditioning. We nevertheless kept the classical Newton method as nonlinear solver. For this simulation, our implementation outperformed MATLAB's PDE toolbox by a factor $4.94$, bringing down our computation time to $78$ s. In this case, the computation time was split more evenly. About $37.7\%$ was used for solving the linear systems while $20.6\%$ was spent in the formation of element matrices.

These experiments reveal the huge computational savings that could be achieved through the implementation of both efficient nonlinear and linear solvers combined with our reassemblage algorithm.

\section{Limitations, potential improvements and future work} \label{se: future_work}
Let us briefly discuss some of the assumptions we have made. 
\begin{itemize}
    \item In Section \ref{se: alternative_assemblage}, we had assumed all elements were of the same type. This allowed us to compute simultaneously the vectorization of all local matrices. Having a mesh with different types of elements is in fact not a limitation. Indeed, a mesh with different types of elements can always be split into several meshes each having a unique type of element. Therefore, the procedure we have described can be applied separately on each single mesh. The different contributions are then merged again before creating the sparse matrix. 
    
    \item In Section \ref{se: classical_assemblage}, we had assumed an affine mapping was used between the reference and the current element. This meant the Jacobian matrix $B_e$ was constant. In practice, other mappings may be used such as isoparametric mappings. In this case, the Jacobian matrix not only depends on the element but also on the quadrature nodes. We can cover for this situation by only changing the definition of the matrix $X$. In particular, one must then compute a matrix of size $n_q \times n_e$ containing all the determinants and another matrix of size $d^2n_q \times n_e$ containing all the vectorizations of the Jacobian matrices.
    
    \item From the very start, in Section \ref{se: fem_basics}, we assumed the coefficient $c$ was a scalar. However, it may happen to be a small matrix of size $d \times d$. This situation arises for instance in case of anisotropic thermal or hydraulic conductivities. Then, returning to point 4 in Section \ref{se: classical_assemblage}, we would have \begin{equation*}
    C_e=\int_{\hat{T}} J_{\hat{\phi}}(\mathbf{\hat{x}})B_e^{-1} \hat{c} B_e^{-T} J_{\hat{\phi}}(\mathbf{\hat{x}})^T |\det(B_e)| \ \mathrm{d} \hat{\Omega} \quad e=1,2,\dots,n_e
\end{equation*}
which would lead to expressing the local conductivity matrix as
\begin{equation*}
    C_e = \sum_{q=1}^{n_q} w_q J_q  B_e^{-1} c_q B_e^{-T} J_q^T |\det(B_e)|  \quad e=1,2,\dots,n_e.  
\end{equation*}
Defining $A_{qe}=B_e^{-1} c_q B_e^{-T}$ and $\lambda_q=w_q |\det(B_e)|$ yields
\begin{equation*}
    \vectorization(C_e) = \sum_{q=1}^{n_q} \lambda_q \ (J_q \otimes J_q)\vectorization(A_{qe}) = \sum_{q=1}^{n_q} \lambda_q \ (J_q \otimes J_q)(B_e^{-1} \otimes B_e^{-1})\vectorization(c_q) \quad e=1,2,\dots,n_e.  
\end{equation*}
Hence, the problem is still successfully decoupled. However, we have not yet worked out the implementation details as anisotropic conductivities are not frequently encountered in practice.

    \item For nonlinear problems where the mesh is changing from one iteration to the next, mesh dependent vectors and matrices will clearly have to be recomputed as well. Provided the algorithms for such tasks are efficient, it should only have a moderate impact on performance. 
\end{itemize}
Despite the already tremendous savings, there is still room for improvements. Let us list a few ideas.
\begin{itemize}
    \item The finite element matrices or tensors we have considered so far always have a certain degree of symmetry. The symmetry naturally carries over to the local matrices and tensors. As we are computing their vectorizations, we straightforwardly know that some of the elements of these vectors will be equal. Having expressed their vectorizations as $V=QX$ allows to avoid unnecessary computations by simply eliminating a few rows of the matrix $Q$ that will only generate duplicates when the product $QX$ is computed. Not only does it induce savings in terms of floating point operations but also in terms of storage because the matrix $Q$ we are storing is smaller than the original one. However, in our implementations, we did not bother with such considerations because neither the product $QX$ nor the storage were the bottlenecks.  
    
    \item The formation of the element mass and conductivity matrices is done sequentially in Algorithm 2 to highlight the similarities in the assembly procedure. However, the formation of all element matrices for each global matrix can naturally be done in parallel. Moreover, the formation of element matrices heavily relies on matrix-matrix operations (level 3 BLAS) that are also well-suited for parallel computations. 

    \item As already noted, if one is not careful enough, much of the computation time for nonlinear problems is spent on repeated calls to MATLAB's built-in \textit{sparse} function. In our implementation, we attempted to reduce the number of calls to a minimum by designing linear operators for carrying out matrix-vector and tensor-vector multiplications. However, we eventually had to form the tangent matrix explicitly in order to use sparse direct solvers (which are used when calling the backslash command). Using GMRES instead allowed to avoid entirely calls to the \textit{sparse} function. However, it was not always advantageous to do so. The reason being that our linear operators could not perform a matrix-vector multiplication as fast as MATLAB could do a sparse matrix-vector multiplication. The use of linear operators was only advantageous when the matrices were very large and GMRES converged very rapidly.
\end{itemize}
Moreover, one might reasonably question whether it is worthwhile assembling the multiplication between tensors and vectors instead of first assembling the tensor and then computing the multiplication. There are two reasons behind this choice:
\begin{enumerate}
    \item Although our procedure successfully applied to $\mathcal{M} \circ_2 (\mathbf{u}-\mathbf{u}^n)$, it did not for $\mathcal{K} \circ_2 \mathbf{u}$. We were unsuccessful in decoupling the terms as we had done for all other quantities. On the other hand, our procedure extended straightforwardly to the assemblage of $\mathcal{M}$ and $\mathcal{C}$.
    
    \item It may happen in applications that the coefficients, although dependent on the unknown, are not too complicated functions. In particular, if the function is linear, then the corresponding tensor is constant and therefore it will only be assembled once. Whereas the multiplication between the tensor and the vector will have to be repeatedly reassembled.
\end{enumerate}

In future work, we will attempt to tackle the following points:
\begin{itemize}
    \item Several interpreted programming languages already offer technologies to speed up critical loops. Python's Numba and MATLAB's Just-In-Time compilation and MEX function generation are just some examples. It might be worthwhile assessing in future work how our strategy compares with existing technologies. Benchmarking would especially target nonlinear problems. Indeed, as we have seen, some of the matrices needed in the assembly process are only computed once, inducing savings in floating point operations.
    
    \item We have already successfully extended our algorithm to the vector PDE of linear elasticity. For the mass matrix, the extension is straightforward due to its Kronecker product structure \cite{voet2020nonlinear}. The stiffness matrix, on the other hand, is much more challenging and we plan to present the algorithm is an upcoming publication.
    
    \item Although all the concepts we have seen apply to 3D problems, the implementation is surely more troublesome. Indeed, dealing with 2D problems allowed us to compute the invariant matrices needed in the algorithm very cheaply. For instance, the vector $\mathbf{b}$ could be computed by simply using the formula for the determinant of $2 \times 2$ matrices along with element-wise operations defined in MATLAB. For 3D problems, these computations become more burdensome. It is therefore worthwhile investigating if these matrices can still be computed efficiently. This is especially important for linear problems. For nonlinear problems, the computation of these vectors and matrices will not play a major role in the overall computational expense since they are usually only computed once.
\end{itemize}

\section{Conclusion} \label{se: conclusion}
In this paper, we have investigated an alternative procedure for assembling finite element matrices and tensors. Although entirely equivalent from a mathematical point of view, our approach allows to design a completely loop-free algorithm. It is clear that the benefits of this approach are especially visible for interpreted programming languages such as MATLAB and Python. But they might hold as well for other languages. Indeed, our approach allows to decouple the various quantities needed in the assemblage of the matrices and therefore avoids having to recompute unnecessary quantities. This strategy becomes particularly appealing for nonlinear problems, when finite element matrices must be reassembled numerous times. Much of the computational workload is concentrated in a few Khatri-Rao products and a small matrix-matrix multiplication. Thanks to its reliance on standard and highly optimized linear algebra operations, our method achieves a tremendous speedup in comparison to traditional assemblage strategies. Finally, at least from a conceptual perspective, our approach is not subjected to any major limitation.  

\section*{Acknowledgments}
\noindent I am grateful to Pablo Antolin for the careful reading of an early version of the manuscript. I would also like to thank two anonymous referees for their constructive comments and helpful suggestions.


\begin{thebibliography}{10}
\expandafter\ifx\csname url\endcsname\relax
  \def\url#1{\texttt{#1}}\fi
\expandafter\ifx\csname urlprefix\endcsname\relax\def\urlprefix{URL }\fi
\expandafter\ifx\csname href\endcsname\relax
  \def\href#1#2{#2} \def\path#1{#1}\fi

\bibitem{ainsworth2011bernstein}
M.~Ainsworth, G.~Andriamaro, O.~Davydov, {Bernstein--B{\'e}zier finite elements
  of arbitrary order and optimal assembly procedures}, SIAM Journal on
  Scientific Computing 33~(6) (2011) 3087--3109.

\bibitem{ainsworth2016bernstein}
M.~Ainsworth, O.~Davydov, L.~L. Schumaker, {Bernstein-B{\'e}zier finite
  elements on tetrahedral--hexahedral--pyramidal partitions}, Computer Methods
  in Applied Mechanics and Engineering 304 (2016) 140--170.

\bibitem{antolin2015efficient}
P.~Antolin, A.~Buffa, F.~Calabro, M.~Martinelli, G.~Sangalli, Efficient matrix
  computation for tensor-product isogeometric analysis: The use of sum
  factorization, Computer Methods in Applied Mechanics and Engineering 285
  (2015) 817--828.

\bibitem{calabro2017fast}
F.~Calabro, G.~Sangalli, M.~Tani, {Fast formation of isogeometric Galerkin
  matrices by weighted quadrature}, Computer Methods in Applied Mechanics and
  Engineering 316 (2017) 606--622.

\bibitem{sangalli2018matrix}
G.~Sangalli, M.~Tani, Matrix-free weighted quadrature for a computationally
  efficient isogeometric k-method, Computer Methods in Applied Mechanics and
  Engineering 338 (2018) 117--133.

\bibitem{mantzaflaris2017low}
A.~Mantzaflaris, B.~J{\"u}ttler, B.~N. Khoromskij, U.~Langer, {Low rank tensor
  methods in Galerkin-based isogeometric analysis}, Computer Methods in Applied
  Mechanics and Engineering 316 (2017) 1062--1085.

\bibitem{scholz2018partial}
F.~Scholz, A.~Mantzaflaris, B.~J{\"u}ttler, {Partial tensor decomposition for
  decoupling isogeometric Galerkin discretizations}, Computer Methods in
  Applied Mechanics and Engineering 336 (2018) 485--506.

\bibitem{hofreither2018black}
C.~Hofreither, A black-box low-rank approximation algorithm for fast matrix
  assembly in isogeometric analysis, Computer Methods in Applied Mechanics and
  Engineering 333 (2018) 311--330.

\bibitem{antolin2020fast}
P.~Antolin, {Fast assembly of Galerkin matrices for 3D solid laminated
  composites using finite element and isogeometric discretizations},
  Computational Mechanics 65~(1) (2020) 135--148.

\bibitem{cuvelier2013efficient}
F.~Cuvelier, C.~Japhet, G.~Scarella, {An efficient way to perform the assembly
  of finite element matrices in Matlab and Octave}, arXiv preprint
  arXiv:1305.3122 (2013).

\bibitem{cuvelier2016efficient}
F.~Cuvelier, C.~Japhet, G.~Scarella, An efficient way to assemble finite
  element matrices in vector languages, BIT Numerical Mathematics 56~(3) (2016)
  833--864.

\bibitem{rahman2013fast}
T.~Rahman, J.~Valdman, {Fast MATLAB assembly of FEM matrices in 2D and 3D:
  Nodal elements}, Applied Mathematics and Computation 219~(13) (2013)
  7151--7158.

\bibitem{anjam2015fast}
I.~Anjam, J.~Valdman, {Fast MATLAB assembly of FEM matrices in 2D and 3D: Edge
  elements}, Applied Mathematics and Computation 267 (2015) 252--263.

\bibitem{quarteroni2009numerical}
A.~Quarteroni, Numerical models for differential problems, Vol.~2, Springer,
  2009.

\bibitem{laursen1978some}
M.~Laursen, M.~Gellert, Some criteria for numerically integrated matrices and
  quadrature formulas for triangles, International journal for numerical
  methods in engineering 12~(1) (1978) 67--76.

\bibitem{hecht2012new}
F.~Hecht, {New development in FreeFem++}, Journal of numerical mathematics
  20~(3-4) (2012) 251--266.

\bibitem{geuzaine2009gmsh}
C.~Geuzaine, J.-F. Remacle, {Gmsh: A 3-D finite element mesh generator with
  built-in pre-and post-processing facilities}, International journal for
  numerical methods in engineering 79~(11) (2009) 1309--1331.

\bibitem{sorber2010khatri_rao}
L.~Sorber,
  \href{https://www.mathworks.com/matlabcentral/fileexchange/28872-khatri-rao-product}{{Khatri-Rao
  product}} (2010).
\newline\urlprefix\url{https://www.mathworks.com/matlabcentral/fileexchange/28872-khatri-rao-product}

\bibitem{ferreira2017validation}
J.~Ferreira, J.~Franssen, T.~Gernay, A.~Gamba, {Validation of
  SAFIR{\textregistered} through DIN EN 1992-1-2 NA} (2017).

\bibitem{saad1986gmres}
Y.~Saad, M.~H. Schultz, {GMRES: A generalized minimal residual algorithm for
  solving nonsymmetric linear systems}, SIAM Journal on scientific and
  statistical computing 7~(3) (1986) 856--869.

\bibitem{voet2020nonlinear}
Y.~Voet, Nonlinear finite elements in dynamics, Tech. rep., {\'E}cole
  polytechnique f{\'e}d{\'e}rale de Lausanne (2020).

\end{thebibliography}
\end{document}